\documentclass[10pt]{article}
\usepackage{amssymb}
\usepackage{amsmath}
\usepackage{amsthm}
\usepackage{hyperref}
\usepackage{cleveref}
\usepackage{bbm} 
\usepackage{thmtools, thm-restate} 
\usepackage{enumitem}   
\usepackage{fullpage} 
\usepackage[colorinlistoftodos,prependcaption]{todonotes}
\usepackage{multirow}
\usepackage{array}

\usepackage{authblk}

\usepackage[scaled=.98,sups]{XCharter} 
\usepackage[charter,scaled=1.07]{newtxmath} 
\usepackage[cal=cm]{mathalfa} 
\DeclareMathAlphabet{\mathsf}{T1}{lmss}{m}{n}
\SetMathAlphabet{\mathsf}{bold}{T1}{lmss}{bx}{n}

\usepackage{titlesec}
\titleformat{\section}
  {\normalfont\scshape\filcenter}{\thesection}{1em}{}
\titleformat{\subsection}
  {\normalsize\scshape}{\thesubsection}{1em}{}

\usepackage{xcolor}
\definecolor{mint}{RGB}{216,254,246}
\definecolor{mintborder}{RGB}{140,220,210}
\definecolor{grayblueborder}{RGB}{100,150,150}
\definecolor{coolgrayborder}{RGB}{120,160,160}
\definecolor{tealborder}{RGB}{0,128,128}

\usepackage{kotex}

\newcommand{\bR}{\mathbb{R}}

\newcommand{\dy}{\mathrm{d} y}
\newcommand{\dz}{\mathrm{d} z}

\newcommand{\gam}{\gamma}

\newcommand{\te}{\theta}

\newcommand{\ome}{\ome}

\newcommand{\Ome}{\Ome}

\renewenvironment{abstract}
  {\begin{center}\begin{minipage}{0.85\linewidth}\vspace{-0.5em}\noindent\textsc{Abstract.} \ignorespaces}
  {\end{minipage}\end{center}}

\numberwithin{equation}{section}

\newtheorem{theorem}{Theorem}
\newtheorem{corollary}{Corollary}
\newtheorem{lemma}{Lemma}
\newtheorem{proposition}{Proposition}
\newtheorem{definition}{Definition}
\newtheorem{remark}{Remark}

\newcommand{\nrm}[1]{\Vert#1\Vert}
\newcommand{\abs}[1]{\left\vert#1\right\vert}
\newcommand{\brk}[1]{\left\langle#1\right\rangle}

\newcommand{\tld}[1]{\widetilde{#1}}

\newcommand{\ud}{\,\mathrm{d}}
\newcommand{\rd}{\partial}

\newcommand{\alp}{\alpha}
\newcommand{\bt}{\beta}
\newcommand{\gma}{\gamma}
\newcommand{\Gma}{\Gamma}
\newcommand{\dlt}{\delta}
\newcommand{\Dlt}{\Delta}

\newcommand{\veps}{\varepsilon}
\newcommand{\kpp}{\kappa}

\newcommand{\tht}{\theta}
\newcommand{\Tht}{\Theta}

\newcommand{\zt}{\zeta}

\newcommand{\vphi}{\varphi}

\newcommand{\bfv}{{\bf v}}

\newcommand{\bbR}{\mathbb R}

\newcommand{\calB}{\mathcal B}

\newcommand{\calE}{\mathcal E}

\newcommand{\calI}{\mathcal I}
\newcommand{\calJ}{\mathcal J}
\newcommand{\calK}{\mathcal K}
\newcommand{\calL}{\mathcal L}
\newcommand{\calM}{\mathcal M}

\newcommand{\calS}{\mathcal S}
\newcommand{\calT}{\mathcal T}

\newcommand{\frkI}{\mathfrak I}

\newcommand\blfootnote[1]{%
  \begingroup
  \renewcommand\thefootnote{}\footnote{#1}%
  \addtocounter{footnote}{-1}%
  \endgroup
}

\newcommand{\BMO}{\mathrm{BMO}}
\newcommand{\MO}{\mathrm{MO}}
\newcommand{\MOP}{{\mathrm{MO}^\vphi}}
\newcommand{\sfd}{\mathsf{d}}

\newcommand{\sfL}{\mathsf{L}}

\newcommand{\dl}{\ud\ell}

\def\XXint#1#2#3{{\setbox0=\hbox{$#1{#2#3}{\int}$}
    \vcenter{\hbox{$#2#3$}}\kern-.5\wd0}}

\usepackage{esint}
\makeatletter
\def\thickhline{%
  \noalign{\ifnum0=`}\fi\hrule \@height \thickarrayrulewidth \futurelet
   \reserved@a\@xthickhline}
\def\@xthickhline{\ifx\reserved@a\thickhline
               \vskip\doublerulesep
               \vskip-\thickarrayrulewidth
             \fi
      \ifnum0=`{\fi}}
\makeatother

\newlength{\thickarrayrulewidth}
\allowdisplaybreaks[4]

\author{\normalsize {Hyunjin In, Dong-ha Kim, Junha Kim, and Seungjae Lee}}
\date{}

\begin{document}

    \title{\large\textbf{A Lagrangian Approach to Uniqueness for Active Scalars in\\
    Generalized Mean Oscillation Spaces}}
    \date\today

    \maketitle

    \blfootnote{
    Hyunjin In: Research Institute of Basic Sciences, Ajou University, 
    Republic of Korea, email: hyunjinin@ajou.ac.kr\\
    Dong-ha Kim(*): Research Institute of Mathematics, Seoul National University,
    Republic of Korea, email: kdh.pde@snu.ac.kr \\
    Junha Kim: Department of Mathematics,
    Ajou University, 
    Republic of Korea, email: junha02@ajou.ac.kr\\
    Seungjae Lee: Department of Mathematical Sciences, Seoul National University, 
    Republic of Korea, email: uhf0125@snu.ac.kr
    \\
    \emph{2020 AMS Mathematics Subject Classification:} 35A02, 76B03, 35Q35, 35Q86, 35Q49, 35Q31 \\
    \emph{Key words: 
    Incompressible Inviscid Fluids, Active Scalar Equations, Uniqueness, Lagrangian viewpoint, Mean Oscillation.} 
    }

    \begin{abstract}
    We develop a unified Lagrangian framework for the uniqueness of weak solutions to a broad class of two-dimensional inviscid active scalar equations in generalized mean oscillation spaces. Under common structural assumptions on the Biot--Savart kernels, we identify mean oscillation spaces that ensure an Osgood condition for the velocity, and establish uniqueness in these classes. For the two-dimensional Euler equations, we establish uniqueness for vorticities whose mean oscillation grows mildly at small scales, extending beyond $\mathrm{BMO}$. For the surface quasi-geostrophic (SQG) and incompressible porous media (IPM) equations, we obtain uniqueness in spaces logarithmically rougher than $C^{0,1}$. We also obtain what are (to our knowledge) the strongest known uniqueness results for the $\alpha$-SQG equation, $0<\alpha<1$, in spaces logarithmically rougher than critical H\"older spaces $C^\alpha$. The framework also applies to various logarithmically perturbed models, yielding uniqueness in the corresponding sharp regularity classes. The proof combines a sharp regularity estimate for the velocity  with a stability estimate for integral operators under measure-preserving rearrangements. This provides a single uniqueness argument for the entire family under common structural assumptions on the kernels.
    \end{abstract}

\section{Introduction} \label{sec1}
\subsection{Active Scalar Equations} \label{sec1.1}
We consider the uniqueness of solutions to the Cauchy problem for a broad class of incompressible inviscid active scalar equations on the whole plane $\bbR^2$:
\begin{equation} \label{eq:active_scalar} \tag{AS} 
    \begin{cases}
        \rd_t\tht + u\cdot\nabla\tht = 0,\\
        u = \calK[\tht],\\
        \tht|_{t=0} = \tht_0.
    \end{cases}
\end{equation}
Here, $\tht: [0,T] \times \bbR^2 \to \bbR$ is the active scalar, $u: [0,T] \times \bbR^2 \to \bbR^2$ is the associated velocity field, and $\calK$ is a Biot--Savart operator. We assume that $\calK$ produces a divergence-free velocity field, $\nabla\cdot u=0$, so that the flow is incompressible. Examples of equations of the form~\eqref{eq:active_scalar} include the two-dimensional Euler equation in vorticity form, the surface quasi-geostrophic (SQG) equation, the $\alp$-surface quasi-geostrophic ($\alp$-SQG) equations for $0<\alp <1$, and the incompressible porous media (IPM) equation, together with logarithmically modified variants of these models.

The uniqueness problem for active scalar equations has been studied from several different perspectives. One classical strategy is Eulerian: given two solutions with the same initial datum, one can compare them by estimating their difference directly in a suitable function space. The foundational result in this direction is Yudovich's uniqueness theorem~\cite{Yudovich_1962} for the two-dimensional Euler equations with bounded vorticity. The proof uses an energy estimate for the velocity difference. This estimate relies on $L^p$-bounds on the velocity gradient and an optimal choice of exponents in Hölder's inequality. Yudovich later extended this mechanism to certain classes of unbounded vorticities characterized by sufficiently slow growth of their $L^p$-norms as $p\to\infty$~\cite{Yudovich_1995Uniqueness}. Vishik~\cite{Vishik_EulerUniqueness} established uniqueness in borderline spaces of Besov type $B_\Gma$, under suitable growth and integrability conditions, using Littlewood--Paley analysis and a paradifferential version of Osgood's argument.

Related Eulerian methods have also been developed for other active scalar equations. For the $\alp$-SQG equations with $0<\alp<1$, Choi, Jung, and Kim~\cite{ChoiJungKim_aSQG} and Jeong~\cite{Jeong_gSQGLectureNote} proved uniqueness in subcritical H\"older spaces $C^{0,\beta}\cap L^1(\bbR^2)$ for $\bt>\alp$ by deriving stability estimates in a negative Sobolev space. For active scalar equations with order-zero velocity operators, including inviscid SQG and IPM, Rusin~\cite{Rusin_Uniqueness} obtained a uniqueness criterion based on BMO control of the scalar gradient under additional assumptions. Azzam and Bedrossian~\cite{AzzamBedrossian_Uniqueness} developed Yudovich-type uniqueness criteria for broad classes of non-diffusive and dissipative active scalar equations, using energy estimates in $L^2$ and $\dot{H}^{-1}$ together with $L^p$--BMO interpolation. More recently, Dom\'inguez and Milman~\cite{DominguezMilman_uniqueness} developed an extrapolation approach to uniqueness for active scalar equations, including SQG and IPM, in Vishik-type spaces $B^1_\Gma$ under an Osgood condition on the associated growth function.

These uniqueness results for the active scalar equations considered here are confined to settings in which the velocity field admits an Osgood modulus of continuity.\footnote{Uniqueness may nevertheless hold in special configurations even when the velocity field does not admit an Osgood modulus of continuity. Such results typically rely on additional geometric or structural information, such as the geometry of the domain and the behavior of the scalar near the boundary~\cite{JeongKimYao_aSQGHalfplane,Zlatos_aSQGHalfplane,AgrawalNahmod_EulerUniquenessRoughDomain}, or the structure of patch-type solutions~\cite{KiselevYaoZlatos_aSQGPatch,CordobaCordobaGancedo_aSQGPatch,Luo_aSQGPatch}, rather than on a general function-space criterion alone.} They suggest that the Osgood condition marks a natural threshold for uniqueness. This perspective is consistent with non-uniqueness results in supercritical classes; for example, the two-dimensional Euler equations~\cite{Vishik_EulerNonuniqueness,Mengual_EulerNonuniqueness,BrueColomboKumar_EulerNonuniqueness}, the inviscid SQG equation~\cite{BuckmasterShkollerVicol_SQGNonuniqueness, Isett_SQGNonuniqueness}, the IPM equation and related active scalar models~\cite{Cordoba_IPMnonuniqueness, IsettVicol_ActiveScalarNonuniqueness}. For further discussion and references, see~\cite{Nonuniquenes_book}.

A complementary approach to \eqref{eq:active_scalar} is provided by the Lagrangian formulation, which exploits the transport structure of the equation. Let $\Phi_t$ denote the flow map associated with velocity field $u$, defined by
\begin{align*}
    \frac{\ud}{\ud t} \Phi_t(x) = u_t(\Phi_t(x)), \qquad \Phi_0(x) = x.
\end{align*}
Here and throughout this paper, time subscripts indicate evaluation at the corresponding time. If $u$ satisfies the Osgood condition, the flow map is uniquely determined. Since $\tht$ is transported by $u$, its value is constant along each flow map:
\begin{equation} \label{eq:Lag_rep}
    \tht_t(x)=\tht_0\bigl(\Phi_t^{-1}(x)\bigr).
\end{equation}
By tracking individual fluid particles, Lagrangian methods have been used to investigate ill-posedness~\cite{JeongKim_EulerIllposedOrlicz, JeongKim_SQGIllposed, BangCheskidov_EulerIllposed, ElgindiJeong_EulerIllposed}, singularity formation~\cite{Elgindi_SingForm, Shkoller_SingForm, Cordoba_SingForm}, and various instability~\cite{ChoiJeong_HillsVortexFilament, DrivasElgindiJeong_Twisting, KiselevSverak_SmallScale, Zlatos_SmallScale}. For the uniqueness problem, this viewpoint leads naturally to a comparison of the flow maps associated with different solutions.

For a passive scalar equation with a prescribed Osgood velocity field, uniqueness follows from the work of Ambrosio and Bernard~\cite{AmbrosioBernard_ODEuniqueness}. Consequently, solutions to such linear transport equations admit a Lagrangian representation. This result also applies to \eqref{eq:active_scalar}: if $u = \calK[\tht]$ satisfies the Osgood condition, then the corresponding solution admits a Lagrangian representation \eqref{eq:Lag_rep}. This alone, however, does not rule out the possible bifurcation of flow maps. Hence, to establish uniqueness for the nonlinear active scalars, one needs to establish that two flow maps associated with any two solutions having the same initial datum coincide.

A classical implementation of this idea for the two-dimensional Euler equations was developed by Marchioro and Pulvirenti~\cite{MarchioroPulvirenti_VortexMethod, MarchioroPulvirenti_Textbook}, whose argument controls a weighted $L^1$--distance between two flow maps. Building on it, Nguyen~\cite{Nguyen_EulerSolutionMap} established continuity of the solution map in the bounded-vorticity class, and Crippa and Stefani~\cite{CrippaStefani_2DYudovich} obtained an elementary Lagrangian proof in uniformly localized Yudovich spaces. Agrawal and Nahmod~\cite{AgrawalNahmod_EulerUniquenessRoughDomain} established uniqueness for bounded vorticity on a class of rough domains by controlling a weighted $L^1$--distance between the flow maps after a novel change of coordinates. This change of coordinates allows them to overcome the loss of log-Lipschitz regularity of the velocity near the boundary. More recently, Drivas and La~\cite{DrivasLa_LagrangianYudovich} established the uniqueness of solutions with vorticities in suitable Yudovich classes by estimating the $L^\infty$--distance between two flow maps and also compared solutions posed on different domains.

Nevertheless, the scope of Lagrangian methods for uniqueness remains less developed than that of Eulerian methods. For the two-dimensional Euler equations, Eulerian uniqueness results extend to broad classes of vorticities in borderline Besov-type spaces, whereas classical Lagrangian arguments have been developed primarily in the Yudovich classes. Moreover, to the best of our knowledge, no Lagrangian proof of uniqueness is available for the $\alp$-SQG, SQG, or IPM equations, even at subcritical regularity.

The purpose of this work is to develop a unified Lagrangian framework for the uniqueness of weak solutions to \eqref{eq:active_scalar} in generalized mean oscillation spaces. Our results cover the two-dimensional Euler equations, the $\alp$-SQG, SQG, and IPM equations, and mild logarithmic perturbations of these models. For each equation, we identify a class within our generalized mean oscillation framework that ensures Osgood continuity of the velocity, and establish uniqueness in this class. The central feature of our approach is that, under common structural assumptions, the dependence on the particular model is encoded in quantitative bounds on its Biot--Savart kernel. This allows us to treat the entire family within a single argument, without requiring separate cancellation properties specific to each equation.

We emphasize that our results concern the uniqueness of weak solutions and do not assert their existence in the same spaces. For example, our results establish uniqueness for the two-dimensional Euler equations with vorticity in $\BMO$, while existence in this class remains a challenging open problem. Our uniqueness classes also encompass $C^\alp$ for the $\alp$-SQG equations and $C^{0,1}$ for the SQG and IPM equations. Nonexistence results are known at some of these regularities; for instance, Choi, Jung, and Kim~\cite{ChoiJungKim_aSQG} proved nonexistence of odd-odd symmetric solutions in $C^\alp$ for the $\alp$-SQG equations, while Jeong and Kim~\cite{JeongKim_SQGIllposed} established nonexistence in $C^{0,1}$ for the SQG equations. 

Nevertheless, uniqueness at low regularity remains meaningful, particularly in nonexistence arguments based on symmetric initial data. For example, since uniqueness in $C^\alp$ was unavailable, \cite{ChoiJungKim_aSQG} established nonexistence only within the class of odd-odd symmetric solutions. Extending this conclusion to arbitrary $C^\alp$ solutions requires uniqueness to ensure that the symmetry of the initial data is preserved. Our uniqueness result for the $\alp$-SQG equations in $C^\alp$ removes this symmetry restriction and implies the nonexistence of any $C^\alp$ solution for the same initial data.

\subsection{Preliminaries} \label{sec1.2}
In this subsection, we introduce the structural assumptions and function spaces used in our main results.

\medskip\noindent\textbf{$\bullet$ Biot--Savart Laws.} \\
Throughout this paper, we consider the active scalar equations with Biot--Savart law $\calK$ satisfying the following structural assumptions.
\begin{enumerate} [label=(A\arabic*)]
    \item \label{A1} For every funciton $f \in C^\infty \cap L^1 (\bbR^2)$, the Biot-Savart law $\calK$ is realized as 
    \begin{align*}
        \calK[f] = \bfv f + K*f,
    \end{align*}
    for some $\bfv \in \bbR^2$ and kernel $K \in C^\infty(\bbR^2\setminus\{0\})$. The convolution with $K$ is interpreted in the principal-value sense when necessary. Also, $\calK[f]$ is divergence-free. 

    \item \label{A2} We define $\kappa : (0,\infty) \to (0,\infty)$ as
        \begin{equation*}
            \kappa(r) := \sup_{|z|=r} |K(z)|.
        \end{equation*}
        Then $\kappa$ satisfies the following properties:
        \begin{enumerate}[label=(\arabic*)]
        \item 
        There exist constants $c_1, c_2 > 0$ satisfying
        \begin{equation*}
            c_1 \kappa(t) \le \kappa(s) \le c_2 \kappa(t), \qquad \text{for all}\quad t/2 \le s \le t < 1.
        \end{equation*}
        \item There exists a constant $C > 0$ such that 
        \begin{equation*}
            \kappa(r) \le C, \qquad \text{for all}\quad r \ge 1.
        \end{equation*}
        \end{enumerate}
    \item \label{A3} The kernel $K$ satisfies the zero-mean condition:
    \begin{equation*}
        \int_{|z|=r} K(z) \,dS(z) = 0, \qquad \text{for all}\quad r > 0.
    \end{equation*}
    \item \label{A4} The kernel $K$ satisfies the following bound for gradient:
    \begin{align*}
        |\nabla K(x)| \leq C \frac{\kappa(|x|)}{|x|},\qquad \text{for all}\quad x \neq 0,
    \end{align*}
    which yields the following difference estimate:
    \begin{equation*}
        |K(x) - K(y)| \le C|x - y| \left( \frac{\kappa(|x|)}{|x|} + \frac{\kappa(|y|)}{|y|} \right), \qquad \text{for all}\quad x \neq y,\quad |x|,|y|>2|x-y|.
    \end{equation*}
\end{enumerate}

We now present several examples of operators $\calK$ satisfying the assumptions above. The most classical example is the Biot--Savart operator for the two-dimensional Euler equations, $\calK_{\text{Euler}}=-\nabla^\perp(-\Dlt)^{-1}$, where the $\perp$ follows the counterclockwise convention. In this case, the associated kernel satisfies the assumptions with $\kpp_{\text{Euler}}(r) \simeq \frac{1}{r}.$ Throughout this paper, $a \simeq b$ denote $\frac{1}{c} a \leq b \leq ca$ for some constant $c>0$. Another example is provided by the SQG equation, for which $\calK_{\text{SQG}}=-\nabla^\perp(-\Dlt)^{-\frac12}$ with the principal-value kernel $K_{\text{SQG}}(x) = \frac{x^\perp}{|x|^3}$. We also consider the incompressible porous media (IPM) equation, whose Biot--Savart law is given by $\calK_{\text{IPM}}[\tht]=-\nabla^\perp \rd_{x_1}(-\Dlt)^{-1}\tht =  \frac{1}{2} e_2 \tht + K_{\text{IPM}}*\tht$, where $K_{\text{IPM}}$ is a principal-value singular kernel given by $K_{\text{IPM}}(x)=\frac{1}{|x|^4} (2x_1x_2, x_2^2-x_1^2)^t$. For both SQG and IPM, the singular kernels correspond to $\kpp(r)\simeq \frac{1}{r^2}$, and the associated convolution operators are understood in the principal-value sense. We further consider the $\alp$-SQG equations, which interpolate between the two-dimensional Euler and SQG equations. For $0<\alp<1$, the corresponding operator is $\calK_{\text{$\alp$-SQG}} = -\nabla^\perp(-\Dlt)^{-1+\frac{\alp}{2}}$, and the associated function $\kpp$ is given by $\kpp_{\text{$\alp$-SQG}}(r) \simeq \frac{1}{r^{1+\alp}}$.

There is also a body of literature on logarithmically perturbed active scalar equations. We do not spell out all the relevant results. A well-studied class of such examples consists of logarithmic perturbations of the Euler equations, with Biot--Savart law $\calK = -m(|D|)\nabla^\perp(-\Dlt)^{-1}$ where $m(|D|)$ is a positive Fourier multiplier exhibiting logarithmic growth or decay. Representative examples can be found in \cite{ChaeConstantin_LogPerturbed,ChaeWu_LogPerturbed, Elgindi_LogPerturbed, Kwon_LogPerturbed}. While these works formulate assumptions on $m(|D|)$ in terms of more general abstract integral conditions, we restrict our attention to concrete examples in which $m(|D|)$ is expressed as a power of a log-function. The corresponding Biot--Savart law also admits a kernel representation, and the associated kernel enjoys a pointwise estimate. For instance, if $m(|D|)=\log^{\gam}\left(1+\frac{1}{|D|}\right)$ for $\gma > 0$, then one readily verifies that the corresponding kernel $K$ has a pointwise bound
\begin{align*}
|K(x)| \leq \frac{C}{|x|}\log^{-\gamma}\left(10+\frac{1}{|x|} \right),\qquad\text{for}\quad 0<|x|<\frac{1}{10}.
\end{align*}
See \cite{Kwon_LogPerturbed} for a related computation. 
Moreover, the kernel $K$ satisfies conditions \ref{A1}--\ref{A4}. Similar logarithmic perturbations can be introduced not only for the Euler equations but also for the other equations discussed above.

\medskip\noindent\textbf{$\bullet$ Generalized Mean Oscillation Spaces.} \\
We consider the following function spaces $\MOP$:
\begin{definition}
    Let $\vphi:(0,10] \to (0,\infty)$ be a given function satisfying the following property:
    there exist constants $c_1, c_2>0$ such that
    \begin{equation} \label{eq:vphi_doubling}
        c_1\vphi(t) \leq \vphi(s) \leq c_2 \vphi(t),
    \end{equation}
    whenever $t/2 \leq s\leq t \leq 10$. We say that a function $\theta \in L^1_{\mathrm{loc}}(\bbR^2)$ belongs to \textit{$\vphi$-mean oscillation space} $\MO^{\vphi}(\bbR^2)$ if there exists $C>0$ such that for $0<r\leq 10$, 
    \begin{equation*}
        \sup_{x\in \bR^2} \fint_{ B(x,r)} \,\abs{\theta(y)-\brk{\theta}_{B(x,r)}}\, \dy 
        \leq C\vphi(r), 
        \qquad \text{where}\,\brk{\theta}_{B(x,r)} :=\fint_{B(x,r)} \theta(y) \ud y.
    \end{equation*}
\end{definition}    
The quantity
\begin{equation*}
    [\te]_{\MO^{\vphi}}:=\sup_{\substack{0<r\leq 10 \\ x \in \bR^2}} \frac{1}{\vphi(r)}\fint_{ B(x,r)} \,\abs{\theta(y)-\brk{\theta}_{B(x,r)}}\,\ud y
\end{equation*}
defines a seminorm on $\MO^{\vphi}$. We remark that $\MO^\vphi$ modulo constants is a Banach space and thus $\MO^\vphi \cap L^p$ for $1\leq p< \infty$ is a Banach space. We emphasize that neither a vanishing nor a boundedness condition is imposed on $\vphi(r)$ as $r\to 0$. \footnote{We remark that in this paper, we only consider $\vphi$ with sub-algebraic growth, that is for each $\veps > 0$,
\[
    \sup_{0<r\leq 10} \vphi(r) r^\veps \leq C_\veps.
\]
This growth restriction ensures the well-definedness of the modulus in \eqref{eq:mu} below. Although this growth assumption may be weakened, we do not pursue its refinement here.}

The most classic example is the space of functions of bounded mean oscillation, $\BMO(\bbR^2)$, corresponding to $\vphi = 1$. A standard borderline example is $f(x)=\log\frac{1}{|x|},$ which belongs to $\BMO(\bbR^2)$ despite being unbounded. In this paper, we also consider spaces allowing mildly unbounded mean oscillation. For $n\geq 2$, define
\begin{align*}
    \log_1(r):=\log\left(10+\frac{1}{r}\right),\qquad\qquad
    \log_{k+1}(r):=\log\left(10+\log_k(r)\right),
    \quad \text{for}\,k\geq 1,
\end{align*}
and set
\begin{align} \label{eq:def_sfL}
    \sfL_n(r):=\prod_{k=2}^{n}\log_k(r).
\end{align}
Note that the product defining $\sfL_n$ starts at $k=2$. A representative choice is $\vphi(r)\sim\sfL_n(r)$ as $r\to 0$. For such $\vphi$, a typical borderline profile near the origin is
\begin{align*}
    f(x)=\log\frac{1}{|x|}\,\sfL_n(|x|).
\end{align*}
More generally, for given $\vphi$ such that $\vphi(r) \to \infty$ as $r \to 0$, a borderline example of $f \in \MOP(\bbR^2)$ is
\begin{align*}
    f(x) = \chi(|x|) \int_{|x|}^1 \frac{\vphi(\ell)}{\ell} \ud \ell,
\end{align*}
where $\chi:[0,\infty) \to \bbR$ is smooth cutoff function supported on $[0,1]$.

We also consider moduli of mean oscillation $\vphi$ that vanish at zero. The corresponding mean oscillation spaces are closely related to spaces defined via a modulus of continuity. More precisely, suppose that $\vphi$ satisfies the Dini condition
\begin{align*}
    \int_0^\veps \frac{\vphi(\ell)}{\ell}\,\ud\ell<\infty,
\end{align*}
for some $\veps>0$. Then every $f\in\MOP$ admits a continuous representative $\bar f$, equal to $f$ almost everywhere in $\bbR^2$, such that
\begin{equation}\label{eq:MO_to_C}
    |\bar f(x)-\bar f(y)| 
    \leq C[f]_{\MO^\vphi}\,\tld{\vphi}(|x-y|),
    \qquad x,y\in\bbR^2,\,0<|x-y|\leq 1,
\end{equation}
where
\begin{equation} \label{eq:tld_vphi}
    \tld{\vphi}(r)
    :=\int_0^r\frac{\vphi(\ell)}{\ell}\,\ud \ell+\vphi(r).
\end{equation}
See \cite[Appendix]{HwangKim_Elliptic} for details.
We henceforth identify $f$ with its continuous representative. Thus, the space $\MO^\vphi$ is continuously embedded into the space $C^{\tld{\vphi}}$, defined by
\begin{align*}
    C^{\tld{\vphi}}(\bbR^2)
    :=\left\{ f\in C^0(\bbR^2) \,\middle|\,
        [f]_{C^{\tld{\vphi}}} := \sup_{0<|x-y|\leq 1}\frac{|f(x)-f(y)|}{\tld{\vphi}(|x-y|)} <\infty
    \right\}.
\end{align*}
We also remark that if $\tld{\vphi} \lesssim \vphi$ then the space $C^\vphi$ and $\MO^\vphi$ coincide. In particular, if $\vphi(r) = r^\alp \sfL_n(r)$ for $0<\alp<1$ and $n \geq 2$ then $\MO^\vphi = C^\vphi$.

\medskip\noindent\textbf{$\bullet$ Uniqueness Classes and Assumptions on $\vphi$.} \\
We now formulate the assumptions on $\vphi$ used in our uniqueness theorem for a fixed Biot--Savart law $\calK$ satisfying \ref{A1}--\ref{A4}. We begin by introducing the modulus associated with the singular-integral part $K*\tht$ of the velocity. Let $\kpp$ be as in \ref{A2}, and define
\begin{equation}\label{eq:mu}
    \mu(r;\kpp,\vphi)
    :=\int_0^r \ell\,\kpp(\ell)\vphi(\ell)\,\ud\ell
      +r\int_r^1\kpp(\ell)\vphi(\ell)\,\ud\ell,
    \qquad r\in(0,1].
\end{equation}
\footnote{One can recover the $\vphi$ for given Osgood modulus $\mu$ and kernel $\kpp$ by $-\mu''(r) = \kpp(r)\vphi(r)$, under suitable assumptions on $\mu$, $\kpp$.} The significance of this modulus is made precise by the following proposition, which is proved in Section~\ref{sec2}.

\begin{proposition} \label{prop:vel_est} 
    Let $K$ satisfies \ref{A2}--\ref{A4} and $\vphi$ satisfies \eqref{eq:vphi_doubling} . For $\tht \in \MOP \cap L^1(\bbR^2)$, we have
    \begin{align*}
        |(K*\tht) (x)- (K*\tht)(y)| \leq C \nrm{\tht}_{\MOP \cap L^1(\bbR^2)} \mu(|x-y|;\kpp,\vphi),
    \end{align*}
    for small enough $|x-y|>0.$
\end{proposition}
From this and \eqref{eq:MO_to_C}, to ensure that $u= \calK[\tht] = \bfv \tht + K*\tht$ satisfies the Osgood condition, it is natural to assume
\begin{enumerate} [label=(\Alph*), ref=(\Alph*), start=2]
    \item \label{B} If $\bfv = 0$, that is $\calK[\tht] = K*\tht$, the induced modulus $\mu(\cdot;\kpp,\vphi)$ in \eqref{eq:mu} satisfies the Osgood condition:
    \begin{equation*}
        \int_0^1 \frac{\ud r}{\mu(r;\kpp,\vphi)} = \infty.
    \end{equation*}
    In addition, if $\bfv \neq 0$, we additionally require that the induced modulus $\tld{\vphi}$ defined in \eqref{eq:tld_vphi} also satisfies the Osgood condition:
    \begin{align*}
        \int_0^1 \frac{\ud r}{\tld{\vphi}(r)} = \infty.
    \end{align*}
\end{enumerate}
We note that by Proposition \ref{prop:vel_est}, for a given Biot--Savart law $\calK$, the condition \ref{B} ensures that, for each $\tht \in \MOP$, the velocity $u = \calK[\tht]$ satisfies the Osgood condition.

\medskip\noindent\textbf{$\bullet$ Weak Solutions to \eqref{eq:active_scalar}.} \\
We now introduce the notion of a weak solution to \eqref{eq:active_scalar}. Suppose that $\calK$ satisfies \ref{A1}--\ref{A4} and that $\vphi$ satisfies \eqref{eq:vphi_doubling} and \ref{B}. A function $\tht \in L^\infty\bigl([0,T];\MO^\vphi \cap L^1(\bbR^2)\bigr)$ is called a weak solution of \eqref{eq:active_scalar} with initial datum $\tht_0 \in \MO^\vphi \cap L^1(\bbR^2)$ if it satisfies
\begin{align*}
    \int_0^T \int_{\bbR^2} \tht \left(\rd_t \zt + u\cdot\nabla \zt \right) \,\ud x\ud t + \int_{\bbR^2}\tht_0(x)\zt(0,x)\,\ud x =0
\end{align*}
for every $\zt\in C_c^\infty([0,T)\times\bbR^2)$.

For given weak solutions $\tht \in L^\infty ([0, T];\MO^\vphi\cap L^1)$, Proposition \ref{prop:vel_est} ensures that the corresponding velocity fields satisfy the Osgood condition. Hence the linear uniqueness theorem (see \cite{BahouriChemin_ODEuniqueness, AmbrosioBernard_ODEuniqueness}, see also \cite{DiPernaLion_ODE}) ensures that the solution $\tht$ must satisfies the Lagrangian structure \eqref{eq:Lag_rep}, that is, for the flow map $\Phi$ generated by $u$, 
\begin{equation*} 
    \tht_t(x) = \tht_0(\Phi_t^{-1} (x)) \qquad \text{for}\quad (t,x) \in [0,T]\times \bbR^2.
\end{equation*}

\subsection{Main Results} \label{sec1.3}

We now state our main uniqueness theorem.

\begin{theorem} \label{thm:thm}
    Assume that $\calK$ satisfies \ref{A1}--\ref{A4}. Also assume that the modulus $\vphi$ satisfies \eqref{eq:vphi_doubling} and \ref{B} for the given operator $\calK$. Given $\tht_0 \in \MOP \cap L^1 (\bbR^2)$ and $T>0$, there exists at most one weak solution to \eqref{eq:active_scalar} in $L^\infty ([0,T];\MOP \cap L^1 (\bbR^2))$.
\end{theorem}

We now describe the scope of our uniqueness theorem and compare it with earlier results. Before turning to individual equations, we record canonical cases of the Osgood condition. Let $\sfL_n$ be defined as in \eqref{eq:def_sfL}. For each integer $n\geq 2$, the following holds:
\begin{equation}\label{eq:mad}
    \vphi(r)\simeq\frac{\sfL_n(r)}{r\kpp(r)}
    \quad\Longrightarrow\quad
    \mu(r;\kpp,\vphi)\simeq r\log(1/r)\sfL_n(r).
\end{equation}
\footnote{If $\mu(r;\kpp,\vphi) = r\log(1/r)\sfL_n(r)$, then since $-\mu''(r) = \kpp(r)\vphi(r)$ the converse holds.}
In particular, when the left-hand condition in \eqref{eq:mad}
holds, $\mu(\cdot;\kpp,\vphi)$ satisfies the Osgood condition.

\medskip\noindent\textbf{$\bullet$ Two-dimensional Euler equations in vorticity form.} \\
The most fundamental equation for \eqref{eq:active_scalar} is the two-dimensional Euler equations in vorticity form. There is substantial uniqueness theory for the Euler equations. We first clarify our uniqueness criterion and then compare it with the previous results in the literature. 

A classical uniqueness result for the two-dimensional Euler equations concerns vorticities with bounded mean oscillation. Theorem \ref{thm:thm} extends uniqueness to vorticities whose mean oscillation may grow mildly at small scales. In particular, by \eqref{eq:mad} with $\kpp(r)\simeq 1/r$, our theorem yields uniqueness for vorticities in $\MOP$ with $\vphi(r)=\sfL_n(r)$ for any fixed $n\geq 2$.

Yudovich~\cite{Yudovich_1995Uniqueness} proved uniqueness theory for unbounded vorticities in the Yudovich Classes. The Yudovich classes are designed to capture the unbounded behavior in terms of the growth of the $L^p$ norm. More precisely, for a positive nondecreasing function $\Theta:[p_0,\infty)\to(0,\infty)$ with $1\leq p_0<\infty$, define
\begin{align*}
    Y^\Theta(\mathbb{R}^2)
    &:=\left\{
        f\in\bigcap_{1\leq p<\infty}L^p(\mathbb{R}^2) \,\middle|\,
        \|f\|_{Y^\Theta} :=\|f\|_{L^1} +\sup_{p\geq p_0}\frac{\|f\|_{L^p}}{\Theta(p)} < \infty
    \right\}.
\end{align*}
Yudovich proved the uniqueness of solutions belonging $L^\infty_t Y^\Tht_x$, provided that
\begin{align} \tag{$\text{B}_\text{Yud}$} \label{eq:B_Yud}
    \int_c^\infty\frac{\ud p}{p\Theta(p)}=\infty, \qquad c\geq p_0.
\end{align}
Unlike the mean oscillation seminorm, $L^p$ norms do not capture the symmetry near singularities. In general, neither $Y^\Theta$ nor $\BMO$ embeds into the other. Nevertheless, allowing the mean oscillation to grow at small scales yields the embedding
\begin{equation}\label{eq:Y_to_MO}
    Y^\Theta(\bbR^2)
    \hookrightarrow \MO^{\vphi_\Theta}(\bbR^2)\cap L^1(\bbR^2),
    \qquad
    \vphi_\Theta(r):=\inf_{p\geq p_0}r^{-2/p}\Theta(p).
\end{equation}
Condition \eqref{eq:B_Yud} implies the Osgood condition for the induced modulus $\mu(\cdot;\kpp,\vphi_\Theta)$ with $\kpp(r)\simeq 1/r$. Together with the embedding \eqref{eq:Y_to_MO}, this shows that Theorem \ref{thm:thm} recovers Yudovich's uniqueness result; see Lemma \ref{lem:Y_to_MO} for details.

Finally, we consider the uniqueness criterion of Vishik~\cite{Vishik_EulerUniqueness} in Besov-type spaces.  For an nondecreasing function $\Gma:\bbR\to[1,\infty)$ satisfying some structural assumptions, define
\begin{align*}
    B_\Gma(\bbR^2)
    :=\left\{ f\in\calS'(\bbR^2)
    \,\middle|\, \nrm{f}_{B_\Gma} := \sup_{N\geq 0} \frac{1}{\Gma(N)} \sum_{j=-1}^{N} \nrm{\Dlt_j f}_{L^\infty} < \infty
    \right\},
\end{align*}
where $(\Dlt_j)_{j\geq-1}$ is an inhomogeneous Littlewood--Paley decomposition and the supremum is over integers $N\geq0$. Vishik proved uniqueness for solutions whose vorticities belong to $L^\infty(0,T;B_\Gma\cap L^p(\bbR^2))$, for some $1<p<2$, provided that
\begin{align} \tag{$\text{B}_\text{Vis}$} \label{eq:B_Vis}
    \int_1^\infty\frac{\mathrm{d}s}{\Gma(s)}=\infty.
\end{align}
A refinement of the classical embedding $\BMO\hookrightarrow B^0_{\infty,\infty}$ yields the continuous embedding
\begin{align} \label{eq:MO_to_B}
    \MO^\vphi \cap L^1(\bbR^2) \hookrightarrow B_{\Gma_\vphi}, 
    \qquad \Gma_\vphi(\alp):=\begin{cases}
        1, & \alpha\leq-1,\\[1mm]
        \displaystyle
        1+ \frac{1}{\log 2}\int_{2^{-(\alpha+1)}}^1 \frac{\vphi(\ell)}\ell\,\ud \ell, & \alpha>-1.
    \end{cases}
\end{align}
For admissible $\vphi$, we further show that the Osgood condition in \ref{B} for $\mu(\cdot;\kpp,\vphi)$, with $\kpp(r)\simeq 1/r$, implies Vishik's criterion \eqref{eq:B_Vis} for $\Gma_\vphi$. Thus, the mean oscillation classes covered by our Euler uniqueness theorem embed into Besov-type spaces satisfying Vishik's criterion; see Lemma \ref{lem:MO_to_B} for details.

\medskip\noindent\textbf{$\bullet$ $\alp$-SQG equations.}\\
There have been several uniqueness results for the $\alp$-SQG equation in subcritical spaces. For example, \cite{Jeong_gSQGLectureNote, ChoiJungKim_aSQG} established the uniqueness of $\alp$-SQG in subcritical H\"older spaces $C^\bt$ with $\bt> \alp$ by obtaining stability estimates in $H^{-1}$ and $H^{-1+\alp}$ respectively. Nevertheless, to the best of our knowledge, uniqueness at the critical H\"older regularity $C^\alp$ has remained open. Our result establishes uniqueness in $C^\alp$ and even in the logarithmically larger space $C^\vphi$, where $\vphi(r) \sim r^\alp \sfL_n(r)$ with $n \geq 2$.

Also, as mentioned earlier, our uniqueness result in $C^\alp$ ensures the propagation of odd-odd symmetry and thus removes the symmetry restriction in the nonexistence result of \cite{ChoiJungKim_aSQG}, yielding nonexistence of arbitrary $C^\alp$ solutions for the same initial data.

\medskip\noindent\textbf{$\bullet$ SQG equations and IPM equations.}\\
There is a substantial literature on the uniqueness of solutions to the SQG and IPM equations. Most existing results rely on mean-oscillation bounds for the gradient of the scalar, $\nabla\tht$. In contrast, our result yields uniqueness in the space of $\MOP$ when $\varphi(r) \sim r\sfL_n(r)$ for some $n\geq 2$, regarding both cases. We remark that this result is not directly comparable with uniqueness results assuming $\nabla\tht\in\BMO$. Indeed, the spaces $C^\vphi$, with $\vphi(r)=r\log\log(1/r)$, and $\BMO^1$ differ substantially. This is because the mean oscillation of the gradient cannot be captured properly in terms of the mean oscillation of the base function. For example, after a suitable cutoff near the origin, the function $x \mapsto|x|\log\log\frac{1}{|x|}$ belongs to $C^\vphi$ but not to $\BMO^1$.

We leave one further remark concerning the IPM equation. Since the Biot–Savart law for IPM contains a Dirac delta term, we must also verify the second condition in \ref{B}, namely the condition on $\widetilde{\varphi}$. This verification is straightforward in our setting. Indeed, when $\varphi(r)=r\sfL_n(r)$, then $\MO^\vphi$ coincides with $C^\vphi$. Therefore, the second condition in \ref{B} is automatically satisfied.

\medskip\noindent\textbf{$\bullet$ Logarithmically Perturbed Model.}\\
As an example of a perturbed model, we consider the logarithmically regularized Euler equation, i.e., the Biot–Savart law is given by $\calK=-m(|D|)\nabla^\perp(-\Dlt)^{-1}$, with $m(|D|)=\log^{\gam}\left(1+\frac{1}{|D|}\right)$ as mentioned in the previous section. For this model, well-posedness and ill-posedness in $H^1$ have been established for different ranges of $\gam >0$; see \cite{ChaeConstantin_LogPerturbed, Kwon_LogPerturbed}. In view of the embedding $H^1(\bR^2)\hookrightarrow\BMO(\bR^2)$, it is natural to ask whether uniqueness also holds in the larger space $\BMO$. Our approach gives positive answer for this whenever $\gamma>0$, and further extends to $\varphi$-mean oscillation spaces whose mean oscillation function $\varphi$ is allowed to grow mildly. We also refer to the log-log singular Euler model in \cite{Elgindi_LogPerturbed}. Our result yields a sharper uniqueness criterion for this model as well. One advantage of our approach is that this proof requires no equation-specific structure and depends only on the asymptotics near the singularity of the corresponding kernel.

\subsection{Proof of the Main Theorems and Idea of Proofs} \label{sec1.4}

In this subsection, we prove our main theorem under the propositions and introduce our key ideas. 

\medskip\noindent\textbf{$\bullet$ Proof of Main Theorem.} 
\begin{proof} [Proof of Theorem~\ref{thm:thm}]
Let $\tht$ and $\rho$ be two weak solutions to \eqref{eq:active_scalar} belonging to $L^\infty([0,T];\MOP(\bbR^2)\cap L^1(\bbR^2))$ with the same initial datum $\tht_0 = \rho_0 \in \MOP(\bbR^2)\cap L^1(\bbR^2)$. Denote their respective velocity fields by $u = \calK[\tht]$ and $v = \calK[\rho]$. By assumption~\ref{B} and Proposition~\ref{prop:vel_est}, the velocity fields $u$ and $v$ satisfy Osgood continuity and therefore generate unique measure-preserving flow maps $\Phi,\Psi:[0,T]\times\bbR^2\to\bbR^2$. Due to \cite{AmbrosioBernard_ODEuniqueness}, $\tht$ and $\rho$ satisfy the Lagrangian representation \eqref{eq:Lag_rep}
\begin{align*}
    \tht_t(x) = \tht_0(\Phi_t^{-1}(x)), \qquad \rho_t(x) = \rho_0(\Psi_t^{-1}(x)), \qquad (t,x) \in [0,T]\times\bbR^2.
\end{align*}
It therefore suffices to prove that $\Phi_t=\Psi_t$ for all $t\in[0,T]$.

Recall that $\calK[f] = \bfv f + K*f$ for $\bfv \in \bbR^2$. The difference between the two flow maps can be written as
\begin{align}\label{eq:basic}
    \begin{split}
        \Phi_t (x) - \Psi_t (x) 
        & = \int_0^t u_s (\Phi_s (x)) - v_s(\Psi_s(x)) \ud s \\
        & = \int_0^t \bfv(\underbrace{\tht_s(\Phi_s(x)) - \rho_s(\Psi_s(x))}_{=\,\tht_0(x)-\rho_0(x)\,=\, 0}) + (K*\tht_s)(\Phi_s(x)) - (K*\rho_s)(\Psi_s(x)) \ud s\\
        & = \int_0^t (K*\tht_s) (\Phi_s (x)) - (K*\tht_s) (\Psi_s (x)) \ud s 
        + \int_0^t (K*\tht_s) (\Psi_s (x)) - (K*\rho_s) (\Psi_s(x)) \ud s.
    \end{split}
\end{align} 
Define the $L^\infty$--distance between the two flow maps at time $t$ by
\begin{equation*}
    \sfd (t) := \nrm{\Phi_t - \Psi_t}_{L^\infty}.
\end{equation*}
We aim to establish the integral inequality
\begin{align} \label{eq:int_ineq}
    \sfd(t) \leq C\int_0^t 
    \mu(\sfd(s); \kpp, \vphi)\,\ud s, \qquad \forall t \in [0,T],
\end{align}
where the constant $C$ depends on $\|\theta_0\|_{\MOP \cap L^1}, \|\tht\|_{L^\infty_t (\MOP\cap L^1)_x}$, and $\|\rho\|_{L^\infty_t (\MOP\cap L^1)_x}$.

By Proposition~\ref{prop:vel_est}, the first integral in \eqref{eq:basic} satisfies
\begin{align*}
    \left| \int_0^t (K*\tht_s) (\Phi_s (x)) - (K*\tht_s) (\Psi_s (x)) \ud s \right|
    \leq C \nrm{\tht}_{L^\infty_t (\MOP\cap L^1)_x} \int_0^t \mu(\sfd(s);\kpp,\vphi) \ud s.
\end{align*}

To estimate the second integral, we use the following proposition, which is a stability estimate for integral operators under measure-preserving rearrangements.
\begin{proposition} \label{prop:vel_diff_est}
    Let $K$ satisfies \ref{A2}--\ref{A4} and $\vphi$ satisfies \eqref{eq:vphi_doubling} and \ref{B} for $\kpp$ corresponding to $K$. Let $\tht \in \MOP(\bbR^2)\cap L^1(\bbR^2)$, and let $\Phi$ and $\Psi$ be  measure-preserving homeomorphisms of $\bbR^2$ with $\nrm{\Phi-\Psi}_{L^\infty}<1$. Suppose that
    \begin{align*}
        \tht\circ\Phi^{-1},\,\tht\circ\Psi^{-1}
        \in\MOP(\bbR^2)\cap L^1(\bbR^2).
    \end{align*}
    Then, 
    \begin{align*}
        \nrm{K*(\tht\circ \Phi^{-1}) - K*(\tht\circ \Psi^{-1})}_{L^\infty} 
        \leq C (\nrm{\tht\circ\Phi^{-1}}_{\MO^\vphi\cap L^1} + \nrm{\tht\circ\Psi^{-1}}_{\MO^\vphi\cap L^1}) \mu(\nrm{\Phi-\Psi}_{L^\infty};\kpp,\vphi).
    \end{align*}
\end{proposition}
The proof of Proposition \ref{prop:vel_diff_est} will be given in Section~\ref{sec3}. By putting $\tht = \tht_0$, $\Phi = \Phi_s$ and $\Psi = \Psi_s$ for $s>0$, Proposition~\ref{prop:vel_diff_est} gives that
\begin{align*}
    |(K*\tht_s) (\Psi_s (x)) - (K*\rho_s) (\Psi_s (x))| \leq C \left( \nrm{\tht_s}_{\MOP\cap L^1} +\nrm{\rho_s}_{\MOP\cap L^1} \right) \mu(\sfd(s);\kpp,\vphi).
\end{align*}
Integrating in time and combining the two estimates, we obtain \eqref{eq:int_ineq}. Since $\mu(\cdot;\kpp,\vphi)$ satisfies the Osgood condition and $\sfd(0)=0$, Osgood's lemma yields $\sfd(t)=0$ for all $t\in[0,T]$. The Lagrangian representations then imply that $\tht=\rho$.
\end{proof}

\noindent\textbf{$\bullet$ Ideas of Proof.} \\
As the proof of the main theorem shows, the argument rests on Propositions~\ref{prop:vel_est} and~\ref{prop:vel_diff_est}. The former follows from potential-theoretic estimates in generalized mean oscillation spaces. The main difficulty lies in the latter, which controls the difference between two velocities in terms of the distance between two flow maps.

To illustrate the idea behind Proposition~\ref{prop:vel_diff_est}, we consider the $\alp$-SQG operator $\calK=-\nabla^\perp(-\Dlt)^{-1+\alp/2}$ with normalized kernel $K(x) = \frac{x^{\perp}}{|x|^{2+\alpha}}$, where $0<\alp<1$, in the critical H\"older space $C^\alp$. We shall aim to show that
\begin{align*}
    \nrm{K*(\tht\circ\Phi^{-1})-K*(\tht\circ\Psi^{-1})}_{L^\infty}
    \leq C\bigl(
        \nrm{\tht\circ\Phi^{-1}}_{C^\alp\cap L^1}
        +\nrm{\tht\circ\Psi^{-1}}_{C^\alp\cap L^1}
    \bigr)
    \nrm{\Phi-\Psi}_{L^\infty}
    \log\frac{1}{\nrm{\Phi-\Psi}_{L^\infty}}.
\end{align*}
This setting is simpler than the general $\MO^\vphi$ framework
but captures the essential mechanism.

Let $\Phi$ and $\Psi$ be measure-preserving homeomorphisms, and assume that $(\tht\circ\Phi^{-1}),(\tht\circ\Psi^{-1})\in C^\alp\cap L^1$. Set $\sfd:=\nrm{\Phi-\Psi}_{L^\infty}$ and suppose that $0<\sfd<1/4$. Fix $x\in\bbR^2$, let $B_r:=B_r(\Psi(x))$ for $r>0$, and set $A:=B_1\setminus B_{2\sfd}$. We decompose the velocity difference at $\Psi(x)$ as 
\begin{align*}
    \bigl(K*(\tht\circ\Phi^{-1})\bigr)(\Psi(x))
    -(K*(\tht\circ\Psi^{-1}))(\Psi(x))
    &=\int_{B_{2\sfd}}K(\Psi(x)-y)
        \bigl(\tht(\Phi^{-1}(y))-\tht(\Psi^{-1}(y))\bigr)\,\ud y\\
    &\quad+\int_A K(\Psi(x)-y)
        \bigl(\tht(\Phi^{-1}(y))-\tht(\Psi^{-1}(y))\bigr)\,\ud y\\
    &\quad+\int_{\bbR^2\setminus B_1}K(\Psi(x)-y)
        \bigl(\tht(\Phi^{-1}(y))-\tht(\Psi^{-1}(y))\bigr)\,\ud y\\
    &=:I_{\mathrm{loc}}+I_{\mathrm{int}}+I_{\mathrm{ext}}.
\end{align*}
The local term $I_{\mathrm{loc}}$ can be controlled by the $C^\alp$ regularity of $\tht\circ\Phi^{-1}$ and $\tht\circ\Psi^{-1}$, using the mean-zero property of $K$ on balls centered at the origin. The exterior term $I_{\mathrm{ext}}$ can be handled using the kernel difference estimate~\ref{A4} and the $L^1$ bounds.

The main challenge is the intermediate term $I_{\mathrm{int}}$. Since the estimate must reflect the displacement $\sfd$, we use measure preservation to change variables. Subtracting the same value $\tht(x)$ from both integrands gives that
\begin{align*}
    I_{\mathrm{int}}
    &=\int_{\Phi^{-1}A}K(\Psi(x)-\Phi(y))
        \bigl(\tht(y)-\tht(x)\bigr)\,\ud y
      -\int_{\Psi^{-1}A} K(\Psi(x)-\Psi(y))
        \bigl(\tht(y)-\tht(x)\bigr)\,\ud y\\
    &=\int_{\Phi^{-1}A\cap \Psi^{-1}A}
        \bigl[K(\Psi(x)-\Phi(y))-K(\Psi(x)-\Psi(y))\bigr]
        \bigl(\tht(y)-\tht(x)\bigr)\,\ud y\\
    &\quad+\int_{\Phi^{-1}A\setminus \Psi^{-1}A}
        K(\Psi(x)-\Phi(y))
        \bigl(\tht(y)-\tht(x)\bigr)\,\ud y\\
    &\quad-\int_{\Psi^{-1}A\setminus\Phi^{-1}A}
        K(\Psi(x)-\Psi(y))
        \bigl(\tht(y)-\tht(x)\bigr)\,\ud y.
\end{align*}
This change of variables allow us to estimate the first term by using the kernel difference bound \ref{A4}, which yields a bound of order $\sfd\log(1/\sfd)$.

The remaining two terms arise because $\Phi^{-1}A$ and $\Psi^{-1}A$ do not coincide. Their integrands cannot be combined pointwise, so we control the regions on which they are supported. Since the flow map need only be H\"older continuous, an argument based on smooth boundaries is unavailable.\footnote{For $\alp$-SQG with scalar regularity $C^\alp$, the velocity is, in general, only log-Lipschitz. The corresponding flow map need not be differentiable.}

Surprisingly, without requiring any regularity of flow maps, the following \textit{geometric lemma} holds
\begin{align*}
    (A\setminus\Phi\Psi^{-1}A)\cup(A\setminus\Psi\Phi^{-1} A)
    \subset
    (B_1\setminus B_{1-\sfd})
    \cup(B_{3\sfd}\setminus B_{2\sfd}).
\end{align*}
See Lemma~\ref{lem_measure} for the precise statement. Thus, after changing variables, both symmetric-difference terms are supported in thin annuli. In particular, we estimate
\begin{align*}
    & \left| \int_{\Psi^{-1}A\setminus \Phi^{-1}A}
    K(\Psi(x)-\Psi(y))
    \bigl(\tht(y)-\tht(x)\bigr)\,\ud y \right| \\
    & \qquad\qquad\qquad\leq \int_{(B_1\setminus B_{1-\sfd}) \cup(B_{3\sfd}\setminus B_{2\sfd})}
    \left| K(\Psi(x)-y)
    \bigl((\tht\circ\Psi^{-1})(y)-(\tht\circ\Psi^{-1})(\Psi(x))\bigr)\right|\,\ud y \\
    & \qquad\qquad\qquad\leq\nrm{\tht\circ\Psi^{-1}}_{C^\alp} \left[ 
    |B_1\setminus B_{1-\sfd}| \sup_{y \in B_1\setminus B_{1-\sfd}} \frac{1}{|\Psi(x)-y|} +  |B_{3\sfd}\setminus B_{2\sfd}| \sup_{y \in B_{3\sfd}\setminus B_{2\sfd}} \frac{1}{|\Psi(x)-y|}   
    \right] \\
    &\qquad\qquad\qquad\leq C\nrm{\tht\circ\Psi^{-1}}_{C^\alp}\,\sfd, 
\end{align*}
and similarly
\begin{align*}
    \left| \int_{\Phi^{-1}A\setminus \Psi^{-1}A}
        K(\Psi(x)-\Phi(y))
        \bigl(\tht(y)-\tht(x)\bigr)\,\ud y \right|
    \leq C\nrm{\tht\circ\Phi^{-1}}_{C^\alp}\,\sfd.
\end{align*}
Combining these bounds with the estimate on the common region gives the desired bound for $I_{\mathrm{int}}$.

The same geometric difficulty arises in mean oscillation spaces, where the absence of pointwise control makes the symmetric-difference terms more delicate. Consider the two-dimensional Euler case with the transported vorticities in $\BMO$. To exploit the $\BMO$ regularity, we divide the intermediate region into dyadic annuli $A_k := B_{2^{k}\sfd}\setminus B_{2^{k-1}\sfd}$, where $2\leq k<N$ and $N\simeq\log(1/\sfd)$, and subtract an appropriate average on each annulus. For each $k$, the integral over $\Phi^{-1}A_k\cap\Psi^{-1}A_k$ can be estimated by the same method. The difficulty lies in the symmetric-difference regions. In the $C^\alpha$ case, pointwise control of the scalar increments gives the required bound; in $\BMO$, the \textit{geometric lemma} alone is insufficient. In particular,  a scale mismatch between the width and radius of these annuli $\Psi^{-1}A_k \triangle \Phi^{-1}A_k$, one only obtains the bounds for each integral of order $\sfd\log(2^k \sfd/\sfd) = \sfd k \log 2$. Summing over $k$ then gives $\sfd(\log(1/\sfd))^2$, which fails the Osgood condition and cannot close the uniqueness argument.

To overcome this difficulty, we examine the geometry of the symmetric-difference regions. The \textit{geometric lemma} gives the decompositions
\begin{align*}
    \Psi^{-1}A_k \setminus \Phi^{-1}A_k 
    = \Psi^{-1} B_{2^{k}\sfd} \setminus \Phi^{-1} B_{2^{k}\sfd} \,\sqcup\, \Phi^{-1} B_{2^{k-1} \sfd} \setminus \Psi^{-1} B_{2^{k-1} \sfd}, \\
    \Phi^{-1}A_k \setminus \Psi^{-1}A_k 
    = \Phi^{-1} B_{2^{k}\sfd} \setminus \Psi^{-1} B_{2^{k}\sfd} \,\sqcup\, \Psi^{-1} B_{2^{k-1} \sfd} \setminus \Phi^{-1} B_{2^{k-1} \sfd},
\end{align*}
In particular, the outer boundary layer of $\Phi^{-1}A_k\setminus\Psi^{-1}A_k$ is the same set as the inner boundary layer of $\Psi^{-1}A_{k+1}\setminus\Phi^{-1}A_{k+1}$, and the corresponding contributions occur with opposite signs. Clearly the same coincidence holds with $\Phi$ and $\Psi$ interchanged. This gives an additional cancellation structure over adjacent annuli. Each paired contribution is $O(\sfd)$, so summing over the scales gives $O(\sfd\log(1/\sfd))$, the Osgood modulus required for uniqueness. The estimate in the general $\MO^\vphi$ setting uses the same cancellation mechanism.

\subsection*{Organization of the paper}
This paper is organized as follows. In Section \ref{sec2}, we present the proof of Proposition \ref{prop:vel_est}, which provides a velocity estimate in mean oscillation spaces. In Section \ref{sec3}, we establish Lemma \ref{lem_measure}, which is crucial for controlling symmetric differences of two preimages of homeomorphisms. Using this lemma, we prove Proposition~\ref{prop:vel_diff_est}, which provides a stability estimate for integral operators under measure-preserving rearrangements.

\subsection*{Acknowledgement}
We thank Professor In-Jee Jeong for his valuable comments.
Hyunjin In was supported by the National Research Foundation of Korea(NRF) grant funded by the Korea government(MSIT) (RS-2024-00360798 and RS-2026-25566720) and Basic Science Research Program through the  National Research Foundation of Korea(NRF) funded by the Ministry of Education(RS-2021-NR060141). Dong-ha Kim was supported by SSEN-Muju Fellow Research Grant and Research Institute of Mathematics, Seoul National University. Junha Kim was supported by the National Research Foundation of Korea(NRF) grant funded by the Korea government(MSIT) (RS-2024-00360798). 
Seungjae Lee was supported by the National Research Foundation of Korea(NRF) grant
from the Korea government (MSIT) (2022R1C1C1011051, RS-2024-00406821).

\subsection*{Declaration of AI Use}

We declare that no AI tools were used in the preparation of this manuscript except for English language editing.

\section{Proof of Proposition \ref{prop:vel_est}} \label{sec2}
In this section, we prove Proposition \ref{prop:vel_est}, which provides the velocity estimate in generalized mean oscillation spaces.
\begin{lemma}\label{lem:basic}
    Let $g \in \MO^{\vphi}(\bbR^2)$. Let $B(x,r) \subset B(\tld{x}, 2r) \subset\bbR^2$ with $0<2r\leq 10$. Then, we have
    \begin{equation*}
        |\brk{g}_{B(x,r)}-\brk{g}_{B(\tld{x}, 2r)}|
        \leq 2^2 [g]_{\MO^{\varphi}}\varphi(2r),
    \end{equation*}
    where $\brk{g}_{B(x,r)} := \fint_{B(x, r)} g(y) \dy$ denotes the average value of $g$ over the ball $B(x,r)$.
\end{lemma}
\begin{proof}
Since $\brk{g}_{\tilde{x}, 2r}$ is constant and $B (x, r) \subset B (\tilde{x}, 2r)$, we obtain
\begin{align*}
    |\brk{g}_{B(x,r)} - \brk{g}_{B(\tld{x}, 2r)}| 
    & = \left| \fint_{B (x, r)} (g(y) - \brk{g}_{B(\tld{x}, 2r)}) \dy \right| \\
    & \le \frac{|B (\tilde{x}, 2r)|}{|B (x, r)|} \fint_{B (\tilde{x}, 2r)} |g(y) - \brk{g}_{B(\tld{x}, 2r)}| \dy \\
    & \le 2^2 [g]_{\MO^{\varphi}}\varphi(2r),
\end{align*}
where we used $|B (\tilde{x}, 2r)|/|B (x, r)| = 2^2$ and applied the definition of the $\MO^{\varphi}$ semi-norm on the ball $B (\tilde{x}, 2r)$.
\end{proof}

\begin{proof} [Proof of Proposition~\ref{prop:vel_est}]
Let $x,y \in \bR^2$ and let $r:=|x-y|$. 
We assume that $0 < r \ll1$. 
Choose an integer $N>0$ such that
\begin{equation*}
    2^{N-1}\leq \frac{1}{r} < 2^N.
\end{equation*}
In particular, $2^N r \simeq 1$.

To utilize the regularity induced from $\MO^\vphi$, we introduce the auxiliary functions

\begin{equation*}
    \phi(z) := \sum_{k=-\infty}^N \brk{\tht}_{B(x,2^kr)} \bold{1}_{B(x,2^kr) \setminus B(x,2^{k-1}r)}(z),
\end{equation*}
and 
\begin{equation*}
    \psi(z) = \sum_{k=-\infty}^0 \brk{\tht}_{B(y,2^kr)} \bold{1}_{B(y,2^kr) \setminus B(y,2^{k-1}r)}(z) + \sum_{k=1}^{N} \brk{\tht}_{B(x,2^kr)} \bold{1}_{B(y,2^kr) \setminus B(y,2^{k-1}r)}(z).
\end{equation*}
Since $\phi$ and $\psi$ are constant on each annulus centered at $x$ and $y$, respectively, the cancellation property of $K$ \ref{A3} yields
\begin{equation*}
    \int_{\bbR^2} K(x-z)\phi(z)\,\ud z=0,
    \qquad
    \int_{\bbR^2} K(y-z)\psi(z)\,\ud z=0.
\end{equation*}

Accordingly, we decompose the velocity difference as follows:
\begin{align*}
    (K*\tht) (x) - (K*\tht) (y) 
    & = \int_{\bbR^2} K(x-z) \tht (z) \,\ud z 
    - \int_{\bbR^2} K(y-z) \tht (z) \,\ud z \\
    & = \int_{\bbR^2} K(x-z) (\tht (z) -\phi(z)) \,\ud z 
    - \int_{\bbR^2} K(y-z) (\tht (z) - \psi(z))  \,\ud z \\
    &=\int_{B(x,2r)} K(x-z) (\tht (z) -\phi(z)) \,\ud z 
    - \int_{B(x,2r)} K(y-z) (\tht (z) - \psi(z)) \,\ud z\\
    &\quad + \int_{B(x,2^{N}r)\setminus B(x,2r)} K(x-z) (\tht (z) -\phi(z)) \,\ud z 
    - \int_{B(x,2^{N}r)\setminus B(x,2r)} K(y-z) (\tht (z) - \psi(z)) \,\ud z\\
    &\quad +\int_{\bbR^2\setminus B(x,2^{N}r)} K(x-z) (\tht (z) -\phi(z)) \,\ud z  
    - \int_{\bbR^2\setminus B(x,2^{N}r)} K(y-z) (\tht (z) - \psi(z)) \,\ud z\\
    &=: L_1+L_2+I_1+I_2+E_1+E_2,
\end{align*}
where $L_j$, $I_j$, and $E_j$ denote the local, intermediate, and exterior contributions, respectively. 

\medskip\noindent\textbf{$\bullet$ Estimates for the local region.} \\
We first estimate $L_1$. By the definition of $\phi$, we have
\begin{align*}
    |L_1| & =\left| \int_{B(x,2r)} K(x-z) (\tht(z) - \phi(z)) \,\ud z \right| \\
    & = \left| \sum_{k=-\infty}^{1} \int_{B(x,2^kr) \setminus B(x,2^{k-1}r)} K(x-z) (\tht(z) - \brk{\tht}_{B(x,2^kr)}) \,\ud z \right|  \\
    & \leq \sum_{k=-\infty}^{1} \int_{B(x,2^kr) \setminus B(x,2^{k-1}r)} \left|K(x-z)\right| \left|\tht(z) - \brk{\tht}_{B(x,2^kr)}\right| \,\ud z \\
    & \leq \sum_{k=-\infty}^{1} \int_{B(x,2^kr) \setminus B(x,2^{k-1}r)} \kappa(|x-z|) \left|\tht(z) - \brk{\tht}_{B(x,2^kr)}\right| \,\ud z \\
    & \leq C \sum_{k=-\infty}^{1} \kappa(2^{k}r) \int_{B(x,2^kr)} \left|\tht(z) - \brk{\tht}_{B(x,2^kr)}\right| \,\ud z \\
    & \leq C[\tht]_{\MO^{\varphi}} \sum_{k=-\infty}^{1} (2^kr)^2 \kappa(2^{k}r) \varphi(2^k r) \\
    & \leq C[\tht]_{\MO^{\varphi}} \int_0^{2r} \ell \kappa(\ell) \vphi(\ell) \,\ud \ell.
\end{align*}
Here and below, we use the doubling conditions of $\kappa$ \ref{A2} and those of $\varphi$ \eqref{eq:vphi_doubling}.

We next estimate $L_2$. Since $r=|x-y|$, we have $B(y,r)\subset B(x,2r)$. Thus, we decompose $B(x,2r)$ into sum of annuli centered at $y$: 
\begin{align*}
    B(x,2r) &= B(y,r) \cup (B(x,2r)\setminus B(y,r)) = \left(\bigcup_{k=-\infty}^0 B(y,2^kr) \setminus B(y,2^{k-1}r) \right) \cup (B(x,2r)\setminus B(y,r)).
\end{align*}
On $B(y,r)$, the definition of $\psi$ therefore yields
\begin{align*}
    \left| \int_{B(y,r)} K(y-z) (\tht(z) - \psi(z)) \,\ud z \right|
    & \leq \sum_{k=-\infty}^{0} \int_{B(y,2^kr) \setminus B(y,2^{k-1}r)} \left| K(y-z) \right| \left| \tht(z) - \brk{\tht}_{B(y,2^kr)} \right| \,\ud z \\
    & \leq C[\tht]_{\MO^{\varphi}} \int_0^{r} \ell \kappa(\ell) \vphi(\ell) \,\ud \ell.
\end{align*}
It remains to estimate the contribution from $B(x,2r)\setminus B(y,r)$.\footnote{For the two-dimensional Euler kernel, this can be controlled easily by using the shell theorem. This cancellation, however, is not needed in our argument.} Since $B(x,2r)\subset B(y,4r)$, we split this region according to the annuli centered at $y$. By the definition of $\psi$, 
\begin{align*}
    \left| \int_{B(x,2r)\setminus B(y,r)} K(y-z) (\tht(z) - \psi(z)) \,\ud z \right| 
    & \leq \left| \int_{(B(x,2r)\setminus B(y,r)) \cap B(y,2r)} K(y-z) (\tht(z) - \brk{\tht}_{B(x,2r)}) \,\ud z \right| \\ 
    & \quad + \left| \int_{(B(x,2r)\setminus B(y,r)) \cap (B(y,4r) \setminus B(y,2r))} K(y-z) (\tht(z) - \brk{\tht}_{B(x,4r)}) \,\ud z \right| \\
    & \leq C\kappa (2r) \int_{B(x,2r)} \left| \tht(z) - \brk{\tht}_{B(x,2r)} \right| \,\ud z \\
    & \quad + C\kappa(4r) \int_{B(x,4r)} \left| \tht(z) - \brk{\tht}_{B(x,4r)} \right| \,\ud z \\
    & \leq C[\tht]_{\MO^{\varphi}} ((2r)^2 \kappa(2r) \varphi(2r) + (4r)^2 \kappa(4r) \varphi(4r)) \\
    & \leq C[\tht]_{\MO^{\varphi}} \int_0^{4r} \ell \kappa(\ell) \vphi(\ell) \,\ud \ell.
\end{align*}
Combining the above estimates, we conclude that 
\begin{equation*}
    |L_1|+|L_2| \leq C[\tht]_{\MO^{\varphi}} \int_0^r \ell\kappa(\ell)\varphi(\ell) \,\ud\ell.
\end{equation*}

\medskip\noindent\textbf{$\bullet$ Estimates for the intermediate terms.} \\
We next estimate the intermediate terms $I_1+I_2$. We recast $I_1+I_2$ as 
\begin{align*}
    I_1+I_2 
    & = \sum_{k=2}^N \int_{B(x,2^{k}r) \setminus B(x,2^{k-1} r)} K(x-z) (\tht(z) - \phi(z)) \ud z - \sum_{k=2}^N\int_{B(x,2^kr) \setminus B(x,2^{k-1}r)} K(y-z) (\tht(z) - \psi(z)) \,\ud z \\
    & = \sum_{k=2}^N \int_{B(x,2^{k} r) \setminus B(x,2^{k-1} r)} (K(x-z) - K(y-z)) (\tht(z) - \phi(z)) \,\ud z \\
    & \quad -\sum_{k=2}^N \int_{B(x,2^{k} r) \setminus B(x,2^{k-1} r)} K(y-z) (\phi(z) - \psi(z)) \,\ud z \\
    & =: J_1+J_2.
\end{align*}

For the term $J_1$, the assumption \ref{A4} yields
\begin{align*}
    |J_1| &\leq \sum_{k=2}^N \int_{B(x,2^{k} r) \setminus B(x,2^{k-1} r)} \left| K(x-z) - K(y-z) \right| \left| \tht(z) - \brk{\tht}_{B(x,2^{k} r)} \right| \,\ud z \\
    &\leq \sum_{k=2}^N \int_{B(x,2^{k} r) \setminus B(x,2^{k-1} r)} C |x-y| \left(\frac{\kappa(|x-z|)}{|x-z|} + \frac{\kappa(|y-z|)}{|y-z|}\right) \left| \tht(z) - \brk{\tht}_{B(x,2^{k} r)} \right| \,\ud z \\
    & \le C r \sum_{k = 2}^N \frac{\kappa(2^{k} r)}{2^{k} r} \int_{B(x,2^{k} r) \setminus B(x,2^{k-1} r)}\left| \theta (z) - \langle \theta \rangle_{B(x,2^kr)} \right| \,\ud z \\
    & \leq C r [\tht]_{\MO^{\varphi}}  \sum_{k=2}^N (2^{k}r) \kappa(2^{k}r) \varphi(2^k r) \\
    & \leq C r [\tht]_{\MO^{\varphi}} \int_{2r}^{2^Nr}\kappa(\ell)\vphi(\ell)\,\ud\ell.
\end{align*}

We now estimate $J_2$. For each $k\geq2$, an annulus $B(x,2^{k} r) \setminus B(x,2^{k-1} r)$ intersects at most three neighboring dyadic annuli centered at $y$. Accordingly, we decompose $J_2$ as follows:
\begin{align*}
    J_2
    & =-\sum_{k=2}^N \int_{(B(x,2^{k} r) \setminus B(x,2^{k-1} r)) \cap (B(y,2^{k-1} r) \setminus B(y,2^{k-2} r))} K(y-z) (\brk{\tht}_{B(x,2^{k} r)} - \brk{\tht}_{B(x,2^{k-1} r)}) \,\ud z \\
    & -\sum_{k=2}^N \int_{(B(x,2^{k} r) \setminus B(x,2^{k-1} r)) \cap (B(y,2^{k} r) \setminus B(y,2^{k-1} r))} K(y-z) (\brk{\tht}_{B(x,2^{k} r)} - \brk{\tht}_{B(x,2^{k} r)}) \,\ud z \\
    & -\sum_{k=2}^{N-1} \int_{(B(x,2^{k} r) \setminus B(x,2^{k-1} r)) \cap (B(y,2^{k+1} r) \setminus B(y,2^{k} r))} K(y-z) (\brk{\tht}_{B(x,2^{k} r)} - \brk{\tht}_{B(x,2^{k+1} r)}) \,\ud z \\
    & - \int_{B(x,2^{N} r) \setminus B(y,2^{N} r)} K(y-z) \brk{\tht}_{B(x,2^N r)} \,\ud z \\
    & =: J_{2,1}+J_{2,2}+J_{2,3}+J_{2,4}.
\end{align*}
Note that the second integral $J_{2,2}$ vanishes. Moreover, since $|x-y|=r$, a simple geometric argument shows that for each $k\geq2$,
\begin{gather*}
    \left| (B(x,2^{k} r) \setminus B(x,2^{k-1} r)) \cap (B(y,2^{k-1} r) \setminus B(y,2^{k-2} r)) \right| 
    + \left| (B(x,2^{k} r) \setminus B(x,2^{k-1} r)) \cap (B(y,2^{k+1} r) \setminus B(y,2^{k} r)) \right| \\
    \leq C2^k r^2.
\end{gather*}
Lemma~\ref{lem:basic} and \eqref{eq:vphi_doubling} imply
\begin{align*}
    \left| \brk{\tht}_{B(x,2^k r)} - \brk{\tht}_{B(x,2^{k-1} r)} \right| + \left| \brk{\tht}_{B(x,2^{k} r)} - \brk{\tht}_{B(x,2^{k+1} r)} \right| \leq C[\tht]_{\MO^{\varphi}} \varphi (2^{k}r).
\end{align*}
Thus, we have
\begin{align*}
    |J_{2,1}|&=\left|\sum_{k=2}^N \int_{(B(x,2^{k} r) \setminus B(x,2^{k-1} r)) \cap (B(y,2^{k-1} r) \setminus B(y,2^{k-2} r))} K(y-z) (\brk{\tht}_{B(x,2^{k} r)} - \brk{\tht}_{B(x,2^{k-1} r)}) \,\ud z \right|\\
    &\leq \sum_{k=2}^N C \left|\brk{\tht}_{B(x,2^{k} r)} - \brk{\tht}_{B(x,2^{k-1} r)}\right| \cdot \left| (B(x,2^{k} r) \setminus B(x,2^{k-1} r)) \cap (B(y,2^{k-1} r) \setminus B(y,2^{k-2} r)) \right| \cdot \kappa(2^k r) \\
    &\leq \sum_{k=2}^N C[\tht]_{\MO^{\varphi}} 2^kr^2 \kappa(2^kr) \varphi(2^{k}r) \\
    & \leq Cr [\tht]_{\MO^{\varphi}} \int_{2r}^{2^Nr}\kappa(\ell)\vphi(\ell)\,\ud\ell.
\end{align*}
Similarly,
\begin{align*}
    |J_{2,3}|
    & = \left| \sum_{k=2}^{N-1} \int_{(B(x,2^{k} r) \setminus B(x,2^{k-1} r)) \cap (B(y,2^{k+1} r) \setminus B(y,2^{k} r))} K(y-z) (\brk{\tht}_{B(x,2^{k} r)} - \brk{\tht}_{B(x,2^{k+1} r)}) \,\ud z \right| \\
    & \leq Cr [\tht]_{\MO^{\varphi}} \int_{2r}^{2^Nr}\kappa(\ell)\vphi(\ell)\,\ud\ell.
\end{align*}
For the term $J_{2,4}$, we note the following elementary geometric estimate: for any
$x,y\in \bR^2$ and $R>0$, 
\begin{equation}\label{lem:geometric}
    |B(x,R)\triangle B(y,R)| \leq C R |x-y|.    
\end{equation}
This implies
\begin{align*}
    |J_{2,4}| 
    & = \left| \int_{B(x,2^{N} r) \setminus B(y,2^{N} r)} K(y-z) \brk{\tht}_{B(x,2^N r)} \,\ud z \right| \\
    &\leq C \kappa(2^Nr) \left|B(x,2^Nr) \triangle B(y,2^Nr)\right| |\brk{\tht}_{B(x,2^Nr)}| \\
    &\leq C \|\theta\|_{L^1}r .
\end{align*}
Therefore, combining the above estimates, we deduce that
\begin{align*}
    |J_2| 
    &\leq Cr [\tht]_{\MO^{\varphi}} \int_{r}^{1} \kappa(\ell) \varphi(\ell) \,\ud \ell + C\|\theta\|_{L^1}r .
\end{align*}

Hence, the intermediate contribution satisfies
\begin{equation*}
    |I_1 + I_2| 
    \leq |J_1| + |J_2| 
    \leq C[\tht]_{\MO^{\varphi}}\, r\int_r^1 \kappa(\ell)\varphi(\ell) \,\ud\ell + C\|\theta\|_{L^1}r .
\end{equation*}

\medskip\noindent\textbf{$\bullet$ Estimates for the exterior terms.} \\
Finally, we estimate the exterior parts $E_1+E_2$.
Since $\phi\equiv 0$ on $\bR^2\setminus B(x,2^N r)$, we can rewrite
\begin{align*}
    E_1+E_2 
    & = \int_{B (x, 2^N r)^c}\left(K(x-z)-K(y-z)\right)\theta(z) \dz + \int_{B (x, 2^N r)^c}K(y-z) \psi(z) \,\dz.\\
    &=: K_1+K_2.
\end{align*}
For $K_1$, the assumption \ref{A4} reveals that
\begin{align*}
    |K_1|
    \leq \int_{B (x, 2^N r)^c} C |x-y| \left(\frac{\kappa(|x-z|)}{|x-z|}+\frac{\kappa(|y-z|)}{|y-z|}\right) |\tht (z)| \,\ud z.
\end{align*}
Since $|x-y| \ll 1$ and $2^N r\sim1$, the factors $\frac{\kappa(|x-z|)}{|x-z|}$ and $\frac{\kappa(|y-z|)}{|y-z|}$
are uniformly bounded on the region of integration. Hence,
\begin{equation*}
    |K_1| \leq C\|\theta\|_{L^1}r.
\end{equation*}\footnote{The $L^1$ condition for $\tht$ can be optimized. For example, for $\kpp(r) = \frac{1}{r^{1+\alp}}$, $\tht \in L^\frac{\alp+2}{2}$ is enough.}

Recalling the definition of $\psi$, we have
\begin{align*}
    K_2= \int_{B(y,2^Nr)\setminus B(x,2^{N}r)}K(y-z) \brk{\tht}_{B(x,2^Nr)} \,\dz.
\end{align*}
As in the estimate for $J_{2,4}$, \eqref{lem:geometric} yields
\begin{align*}
   |K_2| \leq C  \|\theta\|_{L^1}r.
\end{align*}

Therefore, the exterior contribution satisfies
\begin{equation*}
    |E_1+E_2| \leq |K_1| + |K_2| \leq C\|\theta\|_{L^1}r.
\end{equation*}
The proof is complete.
\end{proof}

\section{Proof of Proposition \ref{prop:vel_diff_est}} \label{sec3}

In this section, we prove Proposition \ref{prop:vel_diff_est}, which bounds $\nrm{u-v}_{L^\infty}$ in terms of $\nrm{\Phi-\Psi}_{L^\infty}$. We first record a simple observation, which will play a crucial role in the proof.
\begin{lemma}\label{lem_measure}
Let $\Phi$ and $\Psi$ be homeomorphisms of $\bR^2$ satisfying $\|\Phi-\Psi\|_{L^{\infty}} \leq \delta$.
Let \(D'\subset D\) be sets such that
\begin{equation}\label{assump:1346}
|x-y|>\delta
\qquad\text{for all }x\in D'\text{ and }y\in D^c.
\end{equation}
Then, it holds
    \begin{equation*}
        \Phi^{-1}(D')\cup\Psi^{-1}(D') \subset \Phi^{-1}(D)\cap\Psi^{-1}(D).
    \end{equation*} 
In particular, this implies that
        \begin{equation}\label{implies1354}
         \Phi^{-1}(D)\setminus \Psi^{-1}(D) \subset \Phi^{-1}(D\setminus D'). 
    \end{equation}
\end{lemma}
\begin{proof} 
By symmetry between $\Phi$ and $\Psi$, it suffices to prove that
\begin{equation}\label{eq1355tue}
    \Phi^{-1}(D')\subset \Psi^{-1}(D).
\end{equation}
Suppose, to the contrary, that there exists $z \in \Phi^{-1}(D') \setminus \Psi^{-1}(D)$.  
Then we have
$\Phi(z)\in D'$ and $\Psi(z)\in D^c$, and the assumption \eqref{assump:1346} leads to
\[
|\Phi(z)-\Psi(z)|>\delta.
\] 
This contradicts the assumption $\|\Phi-\Psi\|_{L^{\infty}} \leq \delta$. Indeed, it is obvious to deduce \eqref{implies1354} from \eqref{eq1355tue}. 
The proof is finished.
\end{proof}

Lemma~\ref{lem_measure} directly leads to the following corollary.

\begin{corollary}\label{lem:boundary-layer}
Let $\Phi$ and $\Psi$ be homeomorphisms defined of $\bR^2$ satisfying $\nrm{\Phi-\Psi}_{L^\infty}\leq \dlt$.Then, for every $R\geq\delta$, we have 
\[
    \Phi^{-1}(B_R)\setminus\Psi^{-1}(B_R)
    \subset
    \Phi^{-1}\bigl(B_R\setminus B_{R-\delta}\bigr)
    \cap
    \Psi^{-1}\bigl(B_{R+\delta}\setminus B_R\bigr).
\]
\end{corollary}
\begin{proof}
We use \eqref{implies1354} in Lemma~\ref{lem_measure} with $D=B_R$ and $D' = B_{R-\delta}$ to deduce that
\[
    \Phi^{-1}(B_R)\setminus\Psi^{-1}(B_R)
    \subset
    \Phi^{-1}\bigl(B_R\setminus B_{R-\delta}\bigr).
\]
Applying \eqref{implies1354} again, with the roles of $\Phi$ and
$\Psi$ interchanged, $D'=B_{R+\delta}^c$ and $D=B_R^c$ yields
\[
    \Phi^{-1}(B_R)\setminus\Psi^{-1}(B_R)
    =\Psi^{-1}(B_R^c)\setminus\Phi^{-1}(B_R^c)
    \subset
    \Psi^{-1}\bigl(B_{R}^c\setminus B_{R+\delta}^c\bigr) =\Psi^{-1}\bigl(B_{R+\delta}\setminus B_{R}\bigr).
\]
Combining the two inclusions proves the corollary.
\end{proof}

Now, we are ready to prove Proposition~\ref{prop:vel_diff_est}.
\begin{proof} [Proof of Proposition~\ref{prop:vel_diff_est}]
Throughout the proof, we set
\begin{equation*}
    \sfd:=\nrm{\Phi-\Psi}_{L^\infty}.
\end{equation*}
Since $\Psi$ is bijective, without loss of generality we shall show that
\begin{align*}
    |(K*(\tht\circ\Phi^{-1})) (\Psi(x)) - (K*(\tht\circ\Psi^{-1})) (\Psi(x))| 
    \leq C (\nrm{\tht\circ\Phi^{-1}}_{\MOP \cap L^1} + \nrm{\tht\circ\Psi^{-1}}_{\MOP \cap L^1}) \mu(\sfd;\kpp;\vphi),
\end{align*}
for sufficiently small $\sfd \ll 1$.

For a given $x\in\bbR^2$, we use the shorthand notation
\begin{equation*}
    B_R:=B(\Psi(x),R),
    \qquad
    \brk{\tht\circ\Psi^{-1}}_R
    :=\brk{\tht\circ\Psi^{-1}}_{B_R}
    =\brk{\tht\circ\Psi^{-1}}_{B(\Psi(x),R)},
    \qquad R>0.
\end{equation*}
Since all balls appearing in this proof are centered at $\Psi(x)$, we suppress the center from the notation. We denote dyadic annuli by $$A_k := B_{2^{k}\sfd} \setminus B_{2^{k-1} \sfd}, \qquad k \in \mathbb{Z}.$$

Take $N$ be the integer satisfying
\begin{equation*}
    2^{N-1}\sfd\leq 1<2^N\sfd.
\end{equation*}
That is, $N \simeq \log_2(\frac{1}{\sfd})$. We decompose the difference of the velocity fields as follows:
\begin{align*}
    (K*(\tht\circ\Phi^{-1})) & (\Psi(x)) - (K*(\tht\circ\Psi^{-1})) (\Psi(x)) \\
    & = \int_{\bbR^2} K(\Psi(x)-y) \tht(\Phi^{-1}(y)) \ud y 
    - \int_{\bbR^2} K(\Psi(x)-y) \tht(\Psi^{-1}(y)) \ud y \\
    & =\int_{B_{2\sfd}} K(\Psi(x)-y) \tht(\Phi^{-1}(y)) \ud y 
    - \int_{B_{2\sfd}} K(\Psi(x)-y) \tht(\Psi^{-1}(y)) \ud y\\
    & \quad + \int_{B_{2^N\sfd} \setminus B_{2\sfd}} K(\Psi(x)-y) \tht(\Phi^{-1}(y)) \ud y 
    - \int_{B_{2^N\sfd} \setminus B_{2\sfd}} K(\Psi(x)-y) \tht(\Psi^{-1}(y)) \ud y \\
    & \quad +\int_{\bbR^2\setminus B_{2^N\sfd}} K(\Psi(x)-y) \tht(\Phi^{-1}(y)) \ud y  
    - \int_{\bbR^2\setminus B_{2^N\sfd}} K(\Psi(x)-y) \tht(\Psi^{-1}(y)) \ud y\\
    &=:\calL_1 + \calL_2 + \calI_1 + \calI_2 + \calE_1 + \calE_2,
\end{align*}
where $\calL_j$, $\calI_j$, and $\calE_j$ denote the local, intermediate, and exterior part contributions, respectively.

\medskip\noindent\textbf{$\bullet$ Estimates for the local part.} \\
For the local piece over the ball $B_{2 \sfd}$, we use the mean oscillation bound of $\tht\circ\Phi^{-1}$ and $\tht\circ\Psi^{-1}$. For $\calL_1$, by utilizing the mean-zero condition \ref{A3} of the kernel $K$, we may subtract the average of $\tht\circ\Phi^{-1}$ on each dyadic annulus. Hence,
\begin{align*}
    \calL_1
    &= \int_{B_{2 \sfd}} 
    K (\Psi(x) - y) \,\left(\tht\circ\Phi^{-1}\right)(y) \ud y \\
    &= \sum_{k = -\infty}^1 \int_{A_k} 
    K (\Psi(x) - y) \,\left(\tht\circ\Phi^{-1}\right)(y) \ud y \\
    &= \sum_{k = -\infty}^1 
    \int_{A_k} 
    K (\Psi(x) - y) \left( \left(\tht\circ\Phi^{-1}\right)(y) - \brk{\tht\circ\Phi^{-1}}_{2^{k} \sfd} \right) \ud y.
\end{align*}
Thus, we have
\begin{align*}
|\calL_1|
&\le C 
\sum_{k = -\infty}^1 \kappa(2^{k} \sfd) \int_{A_k} \left| \left(\tht\circ\Phi^{-1}\right) (y) - \brk{\tht\circ\Phi^{-1}}_{{2^{k} \sfd}} \right| \ud y \\
&\le C
[\tht\circ\Phi^{-1}]_{\mathrm{MO}^{\varphi}}
\sum_{k = -\infty}^1 
\kappa(2^{k} \sfd) |B_{2^{k} \sfd}| \varphi(2^{k} \sfd) \\
&\le C
[\tht\circ\Phi^{-1}]_{\mathrm{MO}^{\varphi}}
\sum_{k = -\infty}^1 (2^{k} \sfd)^2 \kappa(2^{k} \sfd) \varphi(2^{k} \sfd) \\
&\le C [\tht\circ\Phi^{-1}]_{\mathrm{MO}^{\varphi}} \int_0^{2 \sfd} \ell \, \kappa(\ell) \varphi(\ell) \ud \ell \\
&\le C [\tht\circ\Phi^{-1}]_{\mathrm{MO}^{\varphi}} \int_0^{\sfd} \ell \, \kappa(\ell) \varphi(\ell) \ud \ell,
\end{align*}
where we have used the doubling properties of $\kappa$ and $\varphi$.

The same argument applied to $\calL_2$ gives
\[
|\calL_2|
= \left| \int_{B_{2 \sfd}} 
K (\Psi(x)-y) \left( \tht\circ\Psi^{-1} \right)(y) \ud y \right| 
\le C [\tht\circ\Psi^{-1} ]_{\mathrm{MO}^{\varphi}}
\int_0^{\sfd} \ell \, \kappa(\ell) \varphi(\ell) \ud \ell.
\]

\noindent

\textbf{$\bullet$ Estimates for the intermediate part.}

We now estimate the intermediate contribution $\calI_1+\calI_2$. We use the mean-zero condition \ref{A3} to write
\begin{align*}
    \calI_1+\calI_2=
    & \int_{B_{2^N\sfd} \setminus B_{2\sfd}} K(\Psi(x)-y) \tht(\Phi^{-1}(y)) \ud y 
    - \int_{B_{2^N\sfd} \setminus B_{2\sfd}} K(\Psi(x)-y) \tht(\Psi^{-1}(y)) \ud y \\ 
    & = \sum_{k=2}^{N} \left[\int_{A_k} K(\Psi (x)-y) \left(\tht(\Phi^{-1}(y)) - \brk{\tht\circ\Psi^{-1}}_{2^{k} \sfd}\right) \ud y 
    - \int_{A_k} K(\Psi (x)-y) \left(\tht(\Psi^{-1}(y)) - \brk{\tht\circ\Psi^{-1}}_{2^{k} \sfd}\right) \ud y \right] \\ 
    & = \sum_{k=2}^{N} \left[\int_{\Phi^{-1} A_k} K(\Psi(x)-\Phi(y)) \left(\tht(y) - \brk{\tht\circ\Psi^{-1}}_{2^{k} \sfd}\right) \ud y 
    - \int_{\Psi^{-1} A_k} K(\Psi(x)-\Psi(y)) \left(\tht(y) - \brk{\tht\circ\Psi^{-1}}_{2^{k} \sfd}\right) \ud y \right] \\ 
    & = \sum_{k=2}^{N} \int_{\Phi^{-1} A_k\cap\Psi^{-1} A_k} \Big( K(\Psi(x)-\Phi(y))-K(\Psi(x)-\Psi(y)) \Big) \left(\tht(y) - \brk{\tht\circ\Psi^{-1}}_{2^{k} \sfd}\right) \ud y \\
    & \quad + \sum_{k=2}^{N} \int_{\Phi^{-1} A_k \setminus \Psi^{-1} A_k} K(\Psi(x)-\Phi(y)) \left(\tht(y) - \brk{\tht\circ\Psi^{-1}}_{2^{k} \sfd}\right) \ud y \\
    & \quad - \sum_{k=2}^{N} \int_{\Psi^{-1} A_k \setminus \Phi^{-1} A_k} K(\Psi(x)-\Psi(y)) \left(\tht(y) - \brk{\tht\circ\Psi^{-1}}_{2^{k} \sfd}\right) \ud y \\
    & =: \frkI_1 + \frkI_2 + \frkI_3. 
\end{align*}

First, the kernel difference estimate \ref{A4} yields the following bound for $\frkI_1$
\begin{align*}
    |\frkI_1|
    & \leq C\sfd \sum_{k=2}^{N} \int_{\Phi^{-1} A_k\cap\Psi^{-1} A_k} \left( \frac{{\kappa}(|\Psi(x)-\Phi(y)|)}{|\Psi(x)-\Phi(y)|} +\frac{{\kappa}(|\Psi(x)-\Psi(y)|)}{|\Psi(x)-\Psi(y)|} \right)  \left| \tht(y) - \brk{\tht\circ\Psi^{-1}}_{2^{k} \sfd} \right| \ud y \\
    & \leq C\sfd \sum_{k=2}^{N} \sup_{z \in  A_k}\left(\frac{\kappa\left( |\Psi(x)-z| \right)}{|\Psi(x)-z|}\right)
    \int_{\Phi^{-1} A_k\cap\Psi^{-1} A_k}   \left| \tht(y) - \brk{\tht\circ\Psi^{-1}}_{2^{k} \sfd} \right| \ud y \\
    & \leq C\sfd \sum_{k=2}^{N} \frac{\kappa(2^k\sfd)}{2^k\sfd}
    \int_{ A_k}   \left| \left(\tht\circ\Psi^{-1}\right)(y) - \brk{\tht\circ\Psi^{-1}}_{2^{k} \sfd} \right| \ud y \\
    & \leq C\sfd [\tht \circ\Psi^{-1}]_{\MOP} \sum_{k=2}^{N} (2^{k}\sfd)\kappa(2^{k}\sfd) \varphi(2^{k+1}\sfd)  \\
    & \leq C\sfd [\tht \circ\Psi^{-1}]_{\MOP} \int_{\sfd}^{1} 
    \kappa(\ell) \varphi(\ell) \dl.
\end{align*}

To estimate $\frkI_2+\frkI_3$, we partition the integration regions and exploit delicate cancellations between corresponding contributions. 
To this end, we define
$$ (\Phi^{-1}A_k \setminus \Psi^{-1}A_k)_{out} := \Phi^{-1} B_{2^{k}\sfd} \setminus \Psi^{-1} B_{2^{k}\sfd}, \qquad (\Phi^{-1}A_k \setminus \Psi^{-1}A_k)_{in} := \Psi^{-1} B_{2^{k-1} \sfd} \setminus \Phi^{-1} B_{2^{k-1} \sfd},$$ for $k\geq 2$.
It is clear that these two sets are disjoint. We claim that the following decomposition holds:
\begin{equation}\label{claim1}
\Phi^{-1}A_k \setminus \Psi^{-1}A_k = (\Phi^{-1}A_k \setminus \Psi^{-1}A_k)_{out} \cup (\Phi^{-1}A_k \setminus \Psi^{-1}A_k)_{in}.
\end{equation}
Indeed, the elementary set identity
\begin{equation*}
    (A\setminus A') \setminus (B\setminus B') = ((A\setminus B) \setminus A') \cup ((B' \setminus A') \cap A),
\end{equation*}
yields that
\begin{equation*}\label{E44}
    \Phi^{-1}A_k \setminus \Psi^{-1}A_k 
    = ((\Phi^{-1} B_{2^{k}\sfd} \setminus \Psi^{-1} B_{2^{k}\sfd}) \setminus \Phi^{-1} B_{2^{k-1} \sfd}) 
    \cup ((\Psi^{-1} B_{2^{k-1} \sfd} \setminus \Phi^{-1} B_{2^{k-1} \sfd}) \cap \Phi^{-1} B_{2^{k}\sfd}).
\end{equation*}
From Corollary~\ref{lem:boundary-layer}, we observe that
\begin{align*}
    \Phi^{-1}B_{2^{k}\sfd} \setminus \Psi^{-1} B_{2^{k}\sfd} \subset
    \Phi^{-1}(B_{2^{k}\sfd} \setminus B_{(2^{k}-1)\sfd}) &= \Phi^{-1}B_{2^{k}\sfd} \setminus \Phi^{-1}B_{(2^{k}-1)\sfd}\\
    \Psi^{-1} B_{2^{k-1} \sfd} \setminus \Phi^{-1}B_{2^{k-1} \sfd} \subset \Phi^{-1}(B_{(2^{k-1}+1)\sfd} \setminus B_{2^{k-1} \sfd}) &= \Phi^{-1}B_{(2^{k-1}+1)\sfd} \setminus \Phi^{-1}B_{2^{k-1} \sfd}.
\end{align*}
This leads to
\begin{align*}
    (\Phi^{-1} B_{2^{k}\sfd} \setminus \Psi^{-1} B_{2^{k}\sfd}) \setminus \Phi^{-1} B_{2^{k-1} \sfd} = \Phi^{-1} B_{2^{k}\sfd} \setminus \Psi^{-1} B_{2^{k}\sfd} &= (\Phi^{-1}A_k \setminus \Psi^{-1}A_k)_{out},\\
    (\Psi^{-1} B_{2^{k-1} \sfd} \setminus \Phi^{-1} B_{2^{k-1} \sfd}) \cap \Phi^{-1} B_{2^{k}\sfd} = \Psi^{-1} B_{2^{k-1} \sfd} \setminus \Phi^{-1} B_{2^{k-1} \sfd} &= (\Phi^{-1}A_k \setminus \Psi^{-1}A_k)_{in},
\end{align*}
respectively, which proves \eqref{claim1}. 
The analogous disjoint decomposition for
$\Psi^{-1}A_k\setminus\Phi^{-1}A_k$ follows in the same way:
\begin{equation*}
    \Psi^{-1}A_k \setminus \Phi^{-1}A_k 
    = \underbrace{\Psi^{-1} B_{2^{k}\sfd} \setminus \Phi^{-1} B_{2^{k}\sfd}}_{=:(\Psi^{-1}A_k \setminus \Phi^{-1}A_k)_{out}} \,\cup\, \underbrace{\Phi^{-1} B_{2^{k-1} \sfd} \setminus \Psi^{-1} B_{2^{k-1} \sfd}}_{=:(\Psi^{-1}A_k \setminus \Phi^{-1}A_k)_{in}}.
\end{equation*}
Notice that
\begin{equation}\label{iden}
    (\Phi^{-1} A_k \setminus \Psi^{-1} A_k)_{out} = (\Psi^{-1} A_{k+1} \setminus \Phi^{-1} A_{k+1})_{in} , \qquad  (\Psi^{-1} A_k \setminus \Phi^{-1} A_k)_{out} = (\Phi^{-1} A_{k+1} \setminus \Psi^{-1} A_{k+1})_{in}.
\end{equation}
These identities \eqref{claim1}-\eqref{iden} allow us to pair the outer contribution from the annuli at scale $2^k\sfd$ with the inner contribution from the annuli at scale $2^{k+1}\sfd$ and obtain an additional cancellation. Consequently, by a telescoping argument, we rewrite $\frkI_2+\frkI_3$ as follows:

\begin{align*}
& \frkI_2+\frkI_3\\
& =\int_{(\Phi^{-1} A_2 \setminus \Psi^{-1} A_2)_{in}} K(\Psi(x)-\Phi(y)) \left(\tht(y) - \brk{\tht\circ\Psi^{-1}}_{2^{2} \sfd}\right) \ud y \\ 
    & \quad - \int_{(\Psi^{-1} A_2 \setminus \Phi^{-1} A_2)_{in}} K(\Psi(x)-\Psi(y)) \left(\tht(y) - \brk{\tht\circ\Psi^{-1}}_{2^{2} \sfd}\right) \ud y \\
    & \quad + \sum_{k=2}^{N-1} \int_{(\Psi^{-1} A_{k+1} \setminus \Phi^{-1} A_{k+1})_{in}} \left[ K(\Psi(x)-\Phi(y)) \left(\tht(y) - \brk{\tht\circ\Psi^{-1}}_{2^{k} \sfd}\right) - K(\Psi(x)-\Psi(y)) \left(\tht(y) - \brk{\tht\circ\Psi^{-1}}_{2^{k+1} \sfd}\right) \right] \ud y \\
    & \quad - \sum_{k=2}^{N-1} \int_{(\Psi^{-1} A_k \setminus \Phi^{-1} A_k)_{out}} \left[ K(\Psi(x)-\Psi(y)) \left(\tht(y) - \brk{\tht\circ\Psi^{-1}}_{2^{k} \sfd}\right) - K(\Psi(x)-\Phi(y)) \left(\tht(y) - \brk{\tht\circ\Psi^{-1}}_{2^{k+1} \sfd}\right) \right] \ud y  \\
    & \quad + \int_{(\Phi^{-1} A_{N} \setminus \Psi^{-1} A_{N})_{out}} K(\Psi(x)-\Phi(y)) \left(\tht(y) - \brk{\tht\circ\Psi^{-1}}_{2^{N} \sfd}\right) \ud y \\
    & \quad - \int_{(\Psi^{-1} A_{N} \setminus \Phi^{-1} A_{N})_{out}} K(\Psi(x)-\Psi(y)) \left(\tht(y) - \brk{\tht\circ\Psi^{-1}}_{2^{N} \sfd}\right) \ud y\\
    &=: \sum_{k=1}^6 \calJ_k.
\end{align*}

For the term $\calJ_1$, we note trivial inclusions $$(\Phi^{-1} A_2 \setminus \Psi^{-1} A_2)_{in} \subset \Phi^{-1}A_2, \qquad (\Phi^{-1} A_2 \setminus \Psi^{-1} A_2)_{in} \subset \Psi^{-1}B_{2\sfd}.$$ 
Hence, we have
\begin{align*}
    |\calJ_1|=& \left|\int_{(\Phi^{-1} A_2 \setminus \Psi^{-1} A_2)_{in}} K(\Psi(x)-\Phi(y)) \left(\tht(y) - \brk{\tht\circ\Psi^{-1}}_{2^{2} \sfd}\right) \ud y\right| \\
    & \leq \sup_{y' \in \Phi^{-1}A_2}|K(\Psi(x)-\Phi(y'))|\int_{\Psi^{-1}B_{2\sfd}}  \left|\tht(y) - \brk{\tht\circ\Psi^{-1}}_{2^{2} \sfd}\right|\ud y \\
    & \leq C\kappa(4\sfd) \int_{B_{4\sfd}} \left| \left( \tht\circ\Psi^{-1} \right)(z) - \brk{\tht\circ\Psi^{-1}}_{2^{2} \sfd}\right| \ud z \\
    & \leq C[ \tht\circ\Psi^{-1}]_{\MOP} \kappa(4\sfd)\sfd^2 \varphi(4\sfd)\\
    & \leq  C[\tht\circ\Psi^{-1}]_{\MOP} \int_{2\sfd}^{4\sfd} \ell\,\kappa(\ell)\varphi(\ell)\ud \ell.
\end{align*} 
By an analogous argument, using
$$(\Psi^{-1} A_2 \setminus \Phi^{-1} A_2)_{in} \subset \Psi^{-1}A_2,$$
we obtain
\begin{align*}
    |\calJ_2|&= \left|\int_{(\Psi^{-1} A_2 \setminus \Phi^{-1} A_2)_{in}} K(\Psi(x)-\Psi(y)) \left(\tht(y) - \brk{\tht\circ\Psi^{-1}}_{2^{2} \sfd}\right) \ud y\right| \\
    & \leq \sup_{y' \in \Psi^{-1}A_2}|K(\Psi(x)-\Psi(y'))|\int_{\Psi^{-1}A_2}  \left|\tht(y) - \brk{\tht\circ\Psi^{-1}}_{2^{2} \sfd}\right|\ud y \\
    & \leq C\kappa(4\sfd) \int_{B_{4\sfd}} \left| \left( \tht\circ\Psi^{-1} \right)(z) - \brk{\tht\circ\Psi^{-1}}_{2^{2} \sfd}\right| \ud z \\
    & \leq  C[\tht\circ\Psi^{-1}]_{\MOP} \int_{2\sfd}^{4\sfd} \ell\,\kappa(\ell)\varphi(\ell)\ud \ell.
\end{align*}

We next split $\calJ_3$ into two parts:
\begin{align*}
    \calJ_3=& \sum_{k=2}^{N-1} \int_{(\Psi^{-1} A_{k+1} \setminus \Phi^{-1} A_{k+1})_{in}} \left[ K(\Psi(x)-\Phi(y))  - K(\Psi(x)-\Psi(y)) \right] \left(\tht(y) - \brk{\tht\circ\Psi^{-1}}_{2^{k+1} \sfd}\right)\ud y \\
    & \quad + \sum_{k=2}^{N-1} \int_{(\Phi^{-1} A_k \setminus \Psi^{-1} A_k)_{out}}  K(\Psi(x)-\Phi(y)) \left( \brk{\tht\circ\Psi^{-1}}_{2^{k+1} \sfd} - \brk{\tht\circ\Psi^{-1}}_{2^{k} \sfd} \right) \ud y.
\end{align*}
For the first part, the kernel difference estimate \ref{A4} yields that 
\begin{align*}
    & \left|\sum_{k=2}^{N-1} \int_{(\Psi^{-1} A_{k+1} \setminus \Phi^{-1} A_{k+1})_{in}} \left[ K(\Psi(x)-\Phi(y))  - K(\Psi(x)-\Psi(y)) \right] \left(\tht(y) - \brk{\tht\circ\Psi^{-1}}_{2^{k+1} \sfd}\right)\ud y\right| \\
    & \leq C\sfd\,\sum_{k=2}^{N-1} \int_{(\Psi^{-1} A_{k+1} \setminus \Phi^{-1} A_{k+1})_{in}} \left(\frac{{\kappa}(|\Psi(x)-\Phi(y)|)}{|\Psi(x)-\Phi(y)|} +\frac{{\kappa}(|\Psi(x)-\Psi(y)|)}{|\Psi(x)-\Psi(y)|} \right) \left|\tht(y) - \brk{\tht\circ\Psi^{-1}}_{2^{k+1} \sfd}\right|\ud y \\
    & \leq C\sfd\,\sum_{k=2}^{N-1} \frac{\kappa(2^k \sfd)}{2^k \sfd} \int_{(\Psi^{-1} A_{k+1} \setminus \Phi^{-1} A_{k+1})_{in}} \left|\tht(y) - \brk{\tht\circ\Psi^{-1}}_{2^{k+1} \sfd}\right| \ud y \\
    & \leq C\sfd\,\sum_{k=2}^{N-1} \frac{\kappa(2^k \sfd)}{2^k \sfd} \int_{A_{k+1}} \left| \left(\tht\circ\Psi^{-1} \right)(y) - \brk{\tht\circ\Psi^{-1}}_{2^{k+1} \sfd}\right| \ud y \\
    & \leq C\sfd\,\sum_{k=2}^{N-1} (2^k\sfd)\kappa(2^k\sfd) \varphi(2^{k+2}\sfd) [\tht\circ\Psi^{-1}]_{\MOP} \\
    & \leq C\sfd\, [\tht\circ\Psi^{-1}]_{\MOP} \int_{\sfd}^{1}  \kappa(\ell) \varphi(\ell)\ud \ell.
\end{align*}
For the second part, it follows from Corollary~\ref{lem:boundary-layer} and the measure-preserving property that
\begin{equation*}
    |(\Phi^{-1} A_k \setminus \Psi^{-1} A_k)_{out}| \leq C2^k \sfd^2.
\end{equation*}
Therefore, 
\begin{align*}
    & \left|\sum_{k=2}^{N-1} \int_{(\Phi^{-1} A_k \setminus \Psi^{-1} A_k)_{out}}  K(\Psi(x)-\Phi(y))  \left( \brk{\tht\circ\Psi^{-1}}_{2^{k+1} \sfd} - \brk{\tht\circ\Psi^{-1}}_{2^{k} \sfd} \right) \ud y \right|\\
    & \qquad\qquad\qquad\qquad\qquad\leq C\sum_{k=2}^{N} \kappa(2^k\sfd) (2^k \sfd^2) \varphi(2^{k+2}\sfd) [\tht\circ\Psi^{-1}]_{\MOP} \\
    & \qquad\qquad\qquad\qquad\qquad\leq C \sfd [\tht\circ\Psi^{-1}]_{\MOP} \int_{\sfd}^{1} \kappa(\ell)\varphi(\ell) \ud \ell.
\end{align*}
Combining the two estimates, we obtain
\begin{equation*}
    |\calJ_3|
    \leq
    C \sfd\,[\tht \circ \Psi^{-1}]_{\MOP}
    \int_{\sfd}^{1}
    \kappa(\ell)\varphi(\ell)\ud\ell.
\end{equation*}
The term $\calJ_4$ is estimated analogously. We omit the details.

Finally, we postpone the estimates of $\calJ_5$ and $\calJ_6$ and treat
them together with the exterior part.

\medskip\noindent\textbf{$\bullet$ Estimates for the exterior part.} \\
For notational simplicity, we write
\[
\calB := B_{2^N \sfd},
\qquad
\calB^c := \bR^2 \setminus \calB.
\]
We decompose the exterior contribution as
\begin{align*}
\calE_1 + \calE_2
&= \int_{\calB^c} K(\Psi(x) - y) \, \left(\tht \circ \Phi^{- 1} \right) (y) \ud y 
- \int_{\calB^c} 
K(\Psi(x) - y) \left( \tht \circ \Psi^{- 1} \right) (y) \ud y \\
&= \int_{\Phi^{- 1} \calB^c \cap \Psi^{- 1} \calB^c} 
\left[ K(\Psi (x) - \Phi (y)) 
- K(\Psi (x) - \Psi (y)) \right] \theta (y) \ud y \\
&\quad + \int_{\Phi^{- 1} \calB^c \setminus \Psi^{- 1} \calB^c} 
K(\Psi (x) - \Phi (y)) 
\theta (y) \ud y \\
&\quad - \int_{\Psi^{- 1} \calB^c \setminus \Phi^{- 1} \calB^c} 
K(\Psi (x) - \Psi (y)) \theta (y) \ud y \\
&=: \calM_1 + \calM_2+\calM_3.
\end{align*}
Since $2^N \sfd \simeq 1$, the kernel difference estimate \ref{A4}
and the uniform bound \ref{A2}-(2) allows us to control $\calM_1$ can be bounded by $L^1$-norm of $\theta$:
\[
|\calM_1|
= \left| \int_{\Phi^{- 1} \calB^c \cap \Psi^{- 1} \calB^c} 
\left[ K(\Psi (x) - \Phi (y)) 
- K(\Psi (x) - \Psi (y)) \right] \theta (y) \ud y \right| 
\le C \nrm{\theta \circ \Psi^{- 1}}_{L^1}\,\sfd,
\]
where we used $\nrm{\theta}_{L^1} = \nrm{\theta \circ \Psi^{- 1}}_{L^1}$.
We next estimate $\calM_2$ and $\calM_3$ together with the terms $\calJ_5$ and $\calJ_6$ from the intermediate part:
\begin{align*}
&\calM_2 + \calM_3 + \calJ_5 + \calJ_6 \\
&= \int_{\Phi^{- 1} \calB^c \setminus \Psi^{- 1} \calB^c} 
K(\Psi (x) - \Phi (y)) \theta (y) \ud y 
- \int_{\Psi^{- 1} \calB^c \setminus \Phi^{- 1} \calB^c} K(\Psi (x) - \Psi (y)) \theta (y) \ud y \\
&\quad + \int_{(\Phi^{- 1}  A_{N} \setminus \Psi^{- 1}  A_{N})_{out}} 
K(\Psi (x) - \Phi (y)) 
(\theta (y) - \brk{\theta \circ \Psi^{- 1}}_{2^N \sfd}) \ud y \\
&\quad - \int_{(\Psi^{- 1}  A_{N} \setminus \Phi^{- 1}  A_{N})_{out}} 
K(\Psi (x) - \Psi (y)) 
(\theta (y) - \brk{\theta \circ \Psi^{- 1}}_{2^N \sfd}) \ud y \\
&= \int_{\Psi^{- 1} \calB \setminus \Phi^{- 1} \calB} 
K(\Psi (x) - \Phi (y)) 
\theta (y) \ud y 
- \int_{\Phi^{- 1} \calB \setminus \Psi^{- 1} \calB} 
K(\Psi (x) - \Psi (y)) \theta (y) \ud y \\
&\quad + \int_{\Phi^{- 1} \calB \setminus \Psi^{- 1} \calB} 
K(\Psi (x) - \Phi (y)) 
(\theta (y) - \brk{\theta \circ \Psi^{- 1}}_{2^N \sfd}) \ud y \\
&\quad - \int_{\Psi^{- 1} \calB \setminus \Phi^{- 1} \calB} 
K(\Psi (x) - \Psi (y)) 
(\theta (y) - \brk{\theta \circ \Psi^{- 1}}_{2^N \sfd}) \ud y \\
&= \int_{\Psi^{- 1} \calB \triangle \Phi^{- 1} \calB} 
\left[ K(\Psi (x) - \Phi (y)) - K(\Psi (x) - \Psi (y)) \right] \theta (y) \ud y \\
&\quad + \int_{\Psi^{- 1} \calB \setminus \Phi^{- 1} \calB} 
K(\Psi (x) - \Psi (y)) \brk{\theta \circ \Psi^{- 1}}_{2^N \sfd} \ud y \\
&\quad - \int_{\Phi^{- 1} \calB \setminus \Psi^{- 1} \calB} 
K(\Psi (x) - \Phi (y)) \brk{\theta \circ \Psi^{- 1}}_{2^N \sfd} \ud y \\
&=: \calT_1 + \calT_2 + \calT_3,
\end{align*}
where $A \triangle B := (A \setminus B) \cup (B \setminus A)$ denotes the symmetric difference of $A$ and $B$.

Since $2^N \sfd \simeq 1$, the kernel difference estimate \ref{A4}
together with the uniform bound \ref{A2}-(2) yields
\[
|\calT_1|
\le C \sfd \, \left( \sup_{\ell \geq \frac12}
\frac{\kappa (\ell)}{\ell} \right)
\int_{\Psi^{- 1} \calB \triangle \Phi^{- 1} \calB}  
|\theta (y)| \ud y 
\le C \nrm{\theta \circ \Psi^{- 1}}_{L^1}\,\sfd.
\]

By Corollary \ref{lem:boundary-layer}, we have
\[
\Psi^{-1} \calB \setminus \Phi^{-1} \calB
\subset 
\Psi^{- 1} (B_{2^N \sfd} \setminus B_{(2^N - 1) \sfd})
\cap
\Phi^{-1} (B_{(2^N + 1) \sfd} \setminus B_{2^N \sfd}).
\]
Since $\Psi$ has the measure-preserving property, it follows that
\[
|\Psi^{-1} \calB \setminus \Phi^{-1} \calB|
\le |\Psi^{- 1} (B_{2^N \sfd} \setminus B_{(2^N - 1) \sfd})|
= |B_{2^N \sfd} \setminus B_{(2^N - 1) \sfd}|
\le C \sfd (2^N \sfd).
\]
From \ref{A2}, it follows that
$|K(\Psi(x) - \Psi(y))| \le C \sup_{\ell \ge 1/2} \kappa(\ell)$.
Moreover, the measure preservation gives
$\nrm{\theta \circ \Psi^{- 1}}_{L^1}=\nrm{\theta}_{L^1}$, and hence
\[
\left|
\brk{\theta \circ \Psi^{- 1}}_{2^N \sfd}
\right|
\le \frac{C}
{(2^N \mathsf d)^2} \nrm{\theta \circ \Psi^{- 1}}_{L^1}.
\]
Consequently,
\begin{align*}
|\calT_2|
&\le \left| \brk{\theta \circ \Psi^{- 1}}_{2^N \sfd} \right|
\left| \int_{\Psi^{- 1} \calB \setminus \Phi^{- 1} \calB} 
K(\Psi (x) - \Psi (y)) \ud y \right| \\
&\le \frac{C}
{(2^N \mathsf d)^2}
(2^N \sfd) \sfd \,
\left( \sup_{\ell \ge \frac12} \kappa (\ell) \right) \nrm{\theta \circ \Psi^{- 1}}_{L^1} \\
&\le C \nrm{\theta \circ \Psi^{- 1}}_{L^1}\,\sfd,
\end{align*}
where the last inequality uses $2^N \sfd > 1$.

The term $\calT_3$ is estimated in the same way.
Hence,
\[
|\calT_3|
\le C \nrm{\theta \circ \Psi^{- 1}}_{L^1}\,\sfd.
\]
We conclude that
\begin{align*}
    |\calE_1 + \calE_2 + \calJ_5 + \calJ_6|
    & = |\calM_1 + \calT_1 + \calT_2 + \calT_3| \\
    & \le C \nrm{\theta \circ \Psi^{- 1}}_{L^1} \,\sfd \\
    & \leq C \nrm{\theta \circ \Psi^{- 1}}_{L^1}\,\mu(\sfd;\kpp,\vphi).
\end{align*}

The proof is complete.

\end{proof}

\newpage\appendix
\setcounter{lemma}{0}
\renewcommand{\thelemma}{\thesection\arabic{lemma}}

\section{Appendix} \label{Appendix}
In this appendix, we prove embeddings between $\MO^\vphi$ and various function spaces. 
We first mention the embedding from generalized mean oscillation spaces $\MO^\vphi$ to modulus of continuity space $C^{\tld{\vphi}}.$ For the proof, see \cite[appendix]{HwangKim_Elliptic}.
\begin{lemma} \label{lem:MO_to_C}
    Let $\vphi: (0,10] \to (0,\infty)$ satisfies \eqref{eq:vphi_doubling} and the Dini condition $\int_0^{\veps} \frac{\vphi(\ell)}{\ell} \ud\ell< \infty$ for some $\veps>0$. Then there is a continuous embedding
    \begin{align} \tag{\ref{eq:MO_to_C}}
        \MO^\vphi(\bbR^2) \hookrightarrow C^{\tld{\vphi}}(\bbR^2),
        \qquad \tld{\vphi}(r) := \int_0^r \frac{\vphi(\ell)}{\ell}\,\ud\ell + \vphi (r),\quad \text{for } r\in[0,10].
    \end{align}
    In particular, if $\vphi$ satisfies
    \begin{align}
        \tld{\vphi}(r) \leq C \vphi(r),
    \end{align}
    then the two spaces coincide, $\MO^\vphi = C^\vphi$.
\end{lemma}

We now verify the embedding from the Yudovich classes to the generalized mean oscillation spaces and compare Yudovich's uniqueness criterion with ours. Recall that the Yudovich space associated with $\theta$ is defined as 
\begin{align*}
    Y^\Theta(\mathbb{R}^2)
    &:=\left\{
        f\in\bigcap_{1\leq p<\infty}L^p(\mathbb{R}^2) \,\middle|\,
        \|f\|_{Y^\Theta} :=\|f\|_{L^1} +\sup_{p\geq p_0}\frac{\|f\|_{L^p}}{\Theta(p)} < \infty
    \right\}.
\end{align*}
\begin{lemma} \label{lem:Y_to_MO}
    Let $1\leq p_0<\infty$, and let $\Tht:[p_0,\infty)\to(0,\infty)$ be nondecreasing. Define
    \begin{equation*}
        \vphi_\Tht(r) :=\inf_{p\geq p_0}r^{-2/p}\Tht(p),
        \qquad r>0.
    \end{equation*}
    Then $\vphi_\Tht$ is positive, continuous, and nonincreasing, with the doubling property \eqref{eq:vphi_doubling}. Moreover, there is a continuous embedding
    \begin{equation} \tag{\ref{eq:Y_to_MO}}
        Y^\Tht(\bbR^2) \hookrightarrow \MO^{\vphi_\Tht}(\bbR^2)\cap L^1(\bbR^2).
    \end{equation}

    The modulus associated with the Euler kernel $\calK = -\nabla^\perp(-\Dlt)^{-1}$ is
    \begin{equation*} 
        \mu(r;\kpp,\vphi_\Tht)
        :=\int_0^r \vphi_\Tht(\ell)\,\ud \ell + r\int_r^1 \frac{\vphi_\Tht(\ell)}{\ell}\,\ud \ell,
        \qquad 0<r\leq 1.
    \end{equation*}
    If $\Tht$ satisfies \eqref{eq:B_Yud}, then $\mu(\cdot;\kpp,\vphi_\Tht)$ satisfies Osgood condition \ref{B}
    \begin{equation*}
        \int_0^1 \frac{\ud r}{\mu(r;\kpp,\vphi_\Tht)}=\infty.
    \end{equation*}
\end{lemma}
\begin{proof}
For each $r>0$, since $\Tht$ is nondecreasing, 
\begin{align*}
    0 < \min \{1,r^{-2/p_0}\}\Tht(p_0)
    \leq \vphi_\Tht(r)
    \leq r^{-2/p_0}\Tht(p_0)<\infty,
\end{align*}
which shows the $\vphi_{\Tht}$ is positive.

For $0 < s \leq r$ and every $p \geq p_0$,
\begin{align} \label{eq:Y_to_MO_bdd1}
    r^{-2/p}\Tht(p)
    \leq s^{-2/p}\Tht(p)
    =\left(\frac rs\right)^{2/p}r^{-2/p}\Tht(p)
    \leq\left(\frac rs\right)^{2/p_0}r^{-2/p}\Tht(p).
\end{align}
Taking the infimum over $p\geq p_0$ gives
\begin{equation*}
    \vphi_\Tht(r)
    \leq\vphi_\Tht(s)
    \leq\left(\frac rs\right)^{2/p_0}\vphi_\Tht(r).
\end{equation*}
This proves monotonicity, continuity on $(0,\infty)$ and doubling property \eqref{eq:vphi_doubling}.

Now we check the embedding from $Y^\Tht$ to $\MO^{\vphi_\Tht}$. Let $f\in Y^\Tht(\bbR^2)$. For every ball $B(x,r)$ and every $p\geq p_0$, the triangle inequality and H\"older's inequality give
\begin{align*}
    \fint_{B(x,r)} \bigl|f-\brk{f}_{B(x,r)}\bigr|\,\ud y
    &\leq 2\fint_{B(x,r)}|f|\,\ud y\\
    &\leq 2|B(x,r)|^{-1/p}\nrm{f}_{L^p}\\
    &\leq 2r^{-2/p}\Tht(p) \nrm{f}_{Y^\Tht}.
\end{align*}
Taking the infimum over $p\geq p_0$, and then the supremum over $x\in \bbR^2$ and $0<r\leq 10$, yields
\begin{align*}
    [f]_{\MO^{\vphi_\Tht}}+\nrm{f}_{L^1} \leq 2\nrm{f}_{Y^\Tht}.
\end{align*}
This proves \eqref{eq:Y_to_MO}.

We now verify that Yudovich's uniqueness criterion \eqref{eq:B_Yud} implies the Osgood condition of $\mu(\cdot;\kpp,\vphi_\Tht)$. 
For $0<r \leq1$ and for every $p\geq \max(p_0,4)$, we have
\begin{align*}
    \mu(r;\kpp,\vphi_\Tht)
    & \leq \Tht(p) \left( \int_0^r s^{-2/p}\,\ud s + r \int_r^1s^{-1-2/p}\,\ud s \right) \\
    & =\Tht(p) \left( \frac{r^{1-2/p}}{1-2/p} +\frac p2\,r\bigl(r^{-2/p}-1\bigr) \right) \\
    &\leq pr^{1-2/p}\Tht(p).
\end{align*}
In the last step we used $p \geq 4$, which implies $(1-2/p)^{-1}\leq p/2$.

In particular, $\mu(r;\kpp,\vphi_\Tht)$ is finite for every
$0<r\leq1$ and tends to zero as $r\downarrow0$.
Since $\vphi_\Tht$ is continuous, differentiation on $(0,1)$ gives
\begin{align*}
    \mu'(r;\kpp,\vphi_\Tht)
    =\int_r^1\frac{\vphi_\Tht(s)}s\,\ud s\geq0,
    \qquad
    \mu''(r;\kpp,\vphi_\Tht)=-\frac{\vphi_\Tht(r)}r<0.
\end{align*}
Thus, its extension by $\mu(0;\kpp,\vphi_\Tht)=0$ is a nondecreasing and concave modulus of continuity.

Set $p_1 := \max(p_0,4)$.
For $0<r\leq e^{-p_1}$, take $p=\log(1/r)$. Since $r^{-2/p}=e^2$, the preceding estimate gives
\begin{equation*}
    \mu(r;\kpp,\vphi_\Tht) \leq e^{2}r \log(1/r) \Tht(\log(1/r)).
\end{equation*} 
Therefore, if \eqref{eq:B_Yud} holds, the change of variables $p=\log(1/r)$ yields
\begin{align*}
    \int_0^{e^{-p_1}} \frac{\ud r}{\mu(r;\kpp,\vphi_\Tht)}
    & \geq C \int_0^{e^{-p_1}} \frac{\ud r}{r\log(1/r)\,\Tht(\log(1/r))}\\
    & = C\int_{p_1}^\infty\frac{\ud p}{p\Tht(p)}
    =\infty.
\end{align*}
This completes the proof.
\end{proof}

We now compare the space $\MO^\vphi$ with Vishik's Besov type space $B_\Gma$. Fix a standard inhomogeneous Littlewood--Paley decomposition $(\Delta_j)_{j\geq-1}$ on $\bbR^2$. Let $\Gma: \bbR \to [1,\infty)$ be a locally Lipschitz continuous monotonically nondecreasing function satisfying the following conditions:
\begin{enumerate} [label=(Vis\arabic*)]
    \item \label{AVis1} $\Gma(\alp)=1$ for every $\alp\leq-1$, and
    $\lim_{\alp\to\infty}\Gma(\alp)=\infty$.
    \item \label{AVis2} There is $C_{\Gma,\mathrm{loc}}\geq1$ such that
    \[
        C_{\Gma,\mathrm{loc}}^{-1}\Gma(\bt)
        \leq\Gma(\alp)
        \leq C_{\Gma,\mathrm{loc}}\Gma(\bt),
        \qquad \alp,\bt\geq-1,\quad |\alp-\bt|\leq1.
    \]
    \item \label{AVis3} There is $C_{\Gma,\mathrm{tail}}>0$ such that
    \[
        \int_{\alp}^{\infty}
            2^{-(\xi-\alp)}\Gma(\xi)\,\ud\xi
        \leq C_{\Gma,\mathrm{tail}}\Gma(\alp),
        \qquad\alp\geq-1.
    \]
\end{enumerate}
Given a function $\Gma:\bbR \to [1,\infty)$ satisfying \ref{AVis1}--\ref{AVis3}, define
\begin{align*}
    B_\Gma(\bbR^2)
    :=\left\{ f\in\calS'(\bbR^2)
    \,\middle|\, \nrm{f}_{B_\Gma} := \sup_{N\geq 0} \frac{1}{\Gma(N)} \sum_{j=-1}^{N} \nrm{\Dlt_j f}_{L^\infty} < \infty
    \right\}.
\end{align*}

\begin{lemma} \label{lem:MO_to_B}
Let $\vphi:(0,10]\to(0,\infty)$ be a continuous function satisfying \eqref{eq:vphi_doubling}. We further assume that $\vphi$ satisfies 
\begin{align}
    & \text{$\vphi$ is nonincreasing}, \label{eq:vphi1}\\
    & \sup_{x \ll 1} \frac{\vphi(x)}{\log \frac{1}{x}} \leq C, \label{eq:vphi2} \\
    & \frac{1}{r} \int_0^r \vphi(\ell) \,\ud \ell 
    \leq C \left(1+\int_r^1 \frac{\vphi(\ell)}{\ell} \,\ud \ell\right) ,\qquad r \in (0,1)\label{eq:vphi3}.
\end{align}
Define
\begin{equation*} 
    \Gma_\vphi(\alpha):=
    \begin{cases}
        1, & \alpha\leq-1,\\[1mm]
        \displaystyle
        1+ \frac{1}{\log 2}\int_{2^{-(\alpha+1)}}^1 \frac{\vphi(\ell)}\ell\,\ud \ell, & \alpha>-1.
    \end{cases}
\end{equation*}
Then $\Gma_\vphi:\bbR\to[1,\infty)$ is a locally Lipschitz, nondecreasing function satisfying Vishik's structural conditions \ref{AVis1}--\ref{AVis3}. Furthermore, there is a continuous embedding
\begin{equation} \tag{\ref{eq:MO_to_B}}
    \MO^\vphi \cap L^1(\bbR^2) \hookrightarrow B_{\Gma_\vphi}(\bbR^2).
\end{equation}

For the Euler kernel $\calK=-\nabla^\perp(-\Dlt)^{-1}$, define
\begin{align*}
    \mu(r;\kpp,\vphi)
    := \int_0^r\vphi(s)\,\ud s + r \int_r^1 \frac{\vphi(s)}s\,\ud s,
    \qquad 0<r\leq1.
\end{align*}
Then the Osgood condition \ref{B} for $\mu$ and Vishik's integral condition \eqref{eq:B_Vis} are equivalent:
\begin{equation*}
    \int_0^1\frac{\ud r}{\mu(r;\kpp,\vphi)}=\infty
    \quad\Longleftrightarrow\quad
    \int_1^\infty\frac{\ud\alpha}{\Gma_\vphi(\alpha)}=\infty.
\end{equation*}
\end{lemma}

\begin{remark}
    The canonical examples of $\vphi$ satisfying \eqref{eq:vphi1}--\eqref{eq:vphi3} are $\vphi = \sfL_n$ for $n \geq 2$.
\end{remark}

\begin{proof}
We first verify that the induced modulus $\Gma_\vphi$ satisfies Vishik's structural assumptions \ref{AVis1}-\ref{AVis3}.  For $\alp>-1$, direct differentiation gives
\begin{align*}
    \Gma_\vphi'(\alp) = \frac{1}{\log 2} \frac{\vphi(2^{-(\alp+1)})}{2^{-(\alp+1)}} 2^{-(\alp+1)} \log 2 = \vphi(2^{-(\alp+1)}) >0.
\end{align*}
This shows that $\Gma_\vphi$ is a nondecreasing, locally Lipschitz function. Moreover, since $\vphi$ is nonincreasing by \eqref{eq:vphi1}, we have $\vphi(2^{-(\alp+1)}) \geq \vphi(1)$, and hence
\begin{align*}
    \Gma_\vphi(\alp) = 1+ \int_{-1}^\alp \Gma_\vphi'(\alp') \ud \alp' \geq 1 + (\alp+1) \vphi(1), \qquad \alp \geq -1.
\end{align*} 
It follows that $\Gma_\vphi(\alp)\to\infty$ as $\alp\to\infty$, which proves \ref{AVis1}.

To verify \ref{AVis2}, let $\alp,\bt\geq-1$ satisfy
$|\alp-\bt|\leq1$. Without loss of generality, we assume that $\alp\leq \bt$. Then,
\begin{align*}
    |\log\Gma_\vphi(\bt) - \log\Gma_\vphi(\alp)|
    & = \left| \int_\alp^\bt \frac{\Gma_\vphi'(\alp')}{\Gma_\vphi(\alp')} \,\ud \alp' \right| \\
    & \leq \int_\alp^\bt \frac{\vphi(2^{-(\alp'+1)})}{1 + (\alp'+1) \vphi(1)} \,\ud\alp' \\
    & \leq |\bt-\alp| \sup_{\alp'\geq -1} \frac{\vphi(2^{-(\alp'+1)})}{1 + (\alp'+1) \vphi(1)} \\
    & \leq C \sup_{x\ll 1} \frac{\vphi(x)}{\log \frac{1}{x}}.
\end{align*}
By \eqref{eq:vphi2}, the last quantity is bounded, proving
\ref{AVis2}.

To verify \ref{AVis3}, we use the identity
\begin{align*}
    \Gma_\vphi(\alp) = 1 + \int_{-1}^\alp \vphi(2^{-(\ell+1)})\,\ud\ell, \qquad \alp\geq -1,
\end{align*}
and write
\begin{align*}
    \int_\alp^\infty 2^{-(\xi-\alp)} \Gma_\vphi(\xi) \,\ud\xi 
    & = \int_\alp^\infty 2^{-(\xi-\alp)} \Gma_\vphi(\alp)\,\ud\xi 
    + \int_\alp^\infty 2^{-(\xi-\alp)} (\Gma_\vphi(\xi) - \Gma_\vphi(\alp)) \,\ud\xi \\
    & = \frac{\Gma_\vphi(\alp)}{\log 2} + \int_\alp^\infty 2^{-(\xi-\alp)} \int_\alp^\xi \vphi(2^{-(\ell+1)})\,\ud\ell\ud\xi \\
    & = \frac{\Gma_\vphi(\alp)}{\log 2} + \int_\alp^\infty \vphi(2^{-(\ell+1)}) \int_\ell^\infty 2^{-(\xi-\alp)} \,\ud\xi\ud\ell \\
    & = \frac{\Gma_\vphi(\alp)}{\log 2} + \frac{1}{\log 2} \int_\alp^\infty \vphi(2^{-(\ell+1)}) 2^{-(\ell-\alp)} \,\ud\ell \\
    & = \frac{\Gma_\vphi(\alp)}{\log 2} + \frac{1}{(\log 2)^2 2^{-(\alp+1)}} \int_0^{2^{-(\alp+1)}} \vphi(\ell) \ud \ell.
\end{align*}
The assumption \eqref{eq:vphi3} ensures that
\begin{align*}
    \frac{1}{2^{-(\alp+1)}} \int_0^{2^{-(\alp+1)}} \vphi(\ell)\,\ud\ell 
    \leq C\left(1+ \int_{2^{-(\alp+1)}}^1 \frac{\vphi(\ell)}{\ell}\,\ud\ell\right) \leq C \Gma_\vphi(\alp).
\end{align*}
This proves \ref{AVis3}.

We next prove the embedding \eqref{eq:MO_to_B}. Let $f\in\MO^\vphi \cap L^1(\bbR^2)$. The two low-frequency blocks satisfies $\nrm{\Dlt_{-1}f}_{L^\infty}+\nrm{\Dlt_0 f}_{L^\infty}\leq C \nrm{f}_{L^1}$. For $j \geq 1$, take standard Littlewood-Paley kernel $h_j(z)=2^{2j}h(2^jz)$, where $h\in\calS(\bbR^2)$ with
$\int h=0$. We write
\begin{align*}
    \Dlt_j f 
    = \int_{\bbR^2} h_j(x-z) (f(z)-\brk{f}_{B(x,2^{-j})}) \,\ud z 
    & = \int_{B(x,2^{-j})} h_j(x-z) (f(z)-\brk{f}_{B(x,2^{-j})}) \,\ud z \\
    & \quad + \sum_{k=0}^{j-1} \int_{B(x,2^{-j+k+1})\setminus B(x,2^{-j+k})} h_j(x-z) (f(z)-\brk{f}_{B(x,2^{-j})}) \,\ud z \\
    & \quad + \int_{B(x,1)^c} h_j(x-z) (f(z)-\brk{f}_{B(x,2^{-j})}) \,\ud z \\
    & =: I_1 + I_2 + I_3.
\end{align*}

We estimate $I_1$ as follows:
\begin{align*}
    \left|\int_{B(x,2^{-j})} h_j(x-z) (f(z)-\brk{f}_{B(x,2^{-j})}) \,\ud z\right|
    & \leq C \nrm{h_j}_{L^\infty} 2^{-2j} [f]_{\MO^\vphi} \vphi(2^{-j}) \\
    & \leq C \nrm{h}_{L^\infty} [f]_{\MO^\vphi} \vphi(2^{-j}). 
\end{align*}

We decompose each summand in $I_2$ as follows:
\begin{align*}
    \int_{B(x,2^{-j+k+1})\setminus B(x,2^{-j+k})} h_j(x-z) & (f(z)-\brk{f}_{B(x,2^{-j})}) \,\ud z \\
    & = \int_{B(x,2^{-j+k+1})\setminus B(x,2^{-j+k})} h_j(x-z) (f(z)-\brk{f}_{B(x,2^{-j+k+1})}) \,\ud z \\
    & \quad + \int_{B(x,2^{-j+k+1})\setminus B(x,2^{-j+k})} h_j(x-z) (\brk{f}_{B(x,2^{-j+k+1})} - \brk{f}_{B(x,2^{-j})}) \,\ud z.
\end{align*}
Using the bound $|h(x)|\leq \frac{C}{(1+|x|)^{100}}$, we obtain
\begin{align*}
    & \left| \int_{B(x,2^{-j+k+1})\setminus B(x,2^{-j+k})} h_j(x-z) (f(z)-\brk{f}_{B(x,2^{-j+k+1})}) \,\ud z \right| \\
    & \qquad\leq C\nrm{h_j}_{L^\infty(B(0,2^{-j+k+1})\setminus B(0,2^{-j+k}))} 2^{-2j+2k} [f]_{\MO^\vphi} \vphi(2^{-j+k+1}) \\
    & \qquad\leq C2^{2j} \nrm{h}_{L^\infty(B(0,2^{k+1})\setminus B(0,2^{k}))} 2^{-2j+2k} [f]_{\MO^\vphi} \vphi(2^{-j+k+1}) \\
    & \qquad\leq \frac{C}{(1+2^k)^{50}} [f]_{\MO^\vphi} \vphi(2^{-j+k+1}) \\
    & \qquad\leq \frac{C}{(1+2^k)^{50}} [f]_{\MO^\vphi} \vphi(2^{-j}),
\end{align*}
where the last inequality follows from the fact that $\vphi$ is nonincreasing. Similarly, by using
\begin{align*}
    |\brk{f}_{B(x,2^{-j+k+1})} - \brk{f}_{B(x,2^{-j})}| \leq C [f]_{\MO^\vphi} \sum_{l=0}^k \vphi(2^{-j+l}) \leq C (k+1) [f]_{\MO^\vphi} \vphi(2^{-j}),
\end{align*}
we estimate the second term as follows:
\begin{align*}
    \left| \int_{B(x,2^{-j+k+1})\setminus B(x,2^{-j+k})} h_j(x-z) (\brk{f}_{B(x,2^{-j+k+1})} - \brk{f}_{B(x,2^{-j})}) \,\ud z \right|
    & \leq C\nrm{h}_{L^\infty(B(0,2^{k+1})\setminus B(0,2^k))} 2^{2k} (k+1) [f]_{\MO^\vphi} \vphi(2^{-j}) \\
    & \leq \frac{C}{(1+2^k)^{50}} [f]_{\MO^\vphi} \vphi(2^{-j}).
\end{align*}
Summing over $k$, we conclude that, for every $j\geq1$,
\begin{align*}
    |I_2| \leq C[f]_{\MO^\vphi} \vphi(2^{-j}).
\end{align*}

For $I_3$, we use $|h_j(x)| \leq C\frac{2^{-2j}}{(1+|x|)^{100}}$ for $|x|\geq 1$ and $j \geq 1$. It follows that
\begin{align*}
    |I_3| 
    & \leq C 2^{-2j} \int_{{B(x,1)}^c} \frac{1}{|x-z|^{100}} |f(z)| \,\ud z + C2^{-2j} |\brk{f}_{B(x,2^{-j})}| \int_{{B(x,1)}^c} \frac{1}{|x-z|^{100}} \ud z \\
    & \leq C 2^{-2j} \nrm{f}_{L^1} + C 2^{-2j} |\brk{f}_{B(x,2^{-j})}|\\
    & \leq C (2^{-2j} + 1) \nrm{f}_{L^1},
\end{align*}
and conclude that $|I_3| \leq C \nrm{f}_{L^1}$. 

Combining the preceding estimates, we obtain, for
$f\in\MO^\vphi\cap L^1(\bbR^2)$ and $j\geq1$,
\begin{align*}
    \nrm{\Dlt_j f}_{L^\infty}  
    & \leq C([f]_{\MO^\vphi} \vphi(2^{-j}) + \nrm{f}_{L^1}) \\
    & \leq C \left( [f]_{\MO^\vphi} + \nrm{f}_{L^1} \right) \vphi(2^{-j}),
\end{align*}
which reveals
\begin{align*}
    \sum_{j=-1}^N\nrm{\Delta_jf}_{L^\infty}
    \leq C \nrm{f}_{\MO^\vphi \cap L^1 (\bbR^2)} \left(1+\sum_{j=0}^N\vphi(2^{-j})\right)
    \leq C \nrm{f}_{\MO^\vphi \cap L^1 (\bbR^2)} \Gma_\vphi(N).
\end{align*}
Dividing by $\Gma_\vphi(N)$ and taking the supremum over integers $N\geq0$ proves \eqref{eq:MO_to_B}.

Finally, we prove the equivalence of the two integral conditions. The definition of $\Gma_\vphi$ gives
\begin{align*}
    \Gma_\vphi\!\left(\frac{\log(1/r)}{\log 2}-1\right)
    =1+\frac{1}{\log 2}\int_r^1\frac{\vphi(s)}s\,\ud s
    \simeq \int_r^1\frac{\vphi(s)}s\,\ud s,
    \qquad \text{for}\quad0< r\leq \frac12.
\end{align*}
Consequently, the assumption \eqref{eq:vphi3} implies that
\begin{align*}
    \mu(r;\kpp,\vphi)
    \simeq
    r\,\Gma_\vphi\!\left(\frac{\log(1/r)}{\log 2}-1\right),
    \qquad \text{for}\quad 0< r\leq \frac12.
\end{align*}
For $0<\veps<1/4$, the change of variables
$\alpha=\log(1/r)/\log 2-1$ therefore yields
\begin{align*}
    \int_\veps^{1/4}
        \frac{\ud r}{\mu(r;\kpp,\vphi)}
    & \simeq \int_\veps^{1/4} \frac{\ud r}{r\, \Gma_\vphi\!\left(\log(1/r)/\log 2-1\right)}\\
    & = \log 2\int_1^{\log(1/\veps)/\log 2-1}
        \frac{\ud\alpha}{\Gma_\vphi(\alpha)}.
\end{align*}
Letting $\veps\downarrow0$ and observing that
\begin{align*}
    \int_{1/4}^1\frac{\ud r}{\mu(r;\kpp,\vphi)}
    \leq \frac{3}{4\mu(1/4;\kpp,\vphi)}<\infty,    
\end{align*}
we conclude that
\begin{align*}
    \int_0^1\frac{\ud r}{\mu(r;\kpp,\vphi)}=\infty
    \quad\Longleftrightarrow\quad
    \int_1^\infty
    \frac{\ud\alpha}{\Gma_\vphi(\alpha)}=\infty.
\end{align*}
This concludes the proof.
\end{proof}

\newpage

\bibliographystyle{alpha}
\bibliography{ref}

\end{document}